\documentclass[12pt,a4paper]{amsart}
\usepackage{a4wide}
\usepackage[english]{babel}
\usepackage[T2A]{fontenc}
\usepackage[utf8]{inputenc}
\usepackage{amsfonts}
\usepackage{amssymb, amsthm, amscd}
\usepackage{amsmath}
\usepackage{mathtools}
\usepackage{needspace}
\usepackage{etoolbox}
\usepackage{lipsum}
\usepackage{comment}
\usepackage{cmap}
\usepackage[pdftex]{graphicx}
\usepackage[unicode]{hyperref}
\usepackage[matrix,arrow,curve]{xy}
\usepackage[usenames,dvipsnames]{xcolor}
\usepackage{colortbl}
\usepackage{textcomp}
\usepackage{cite}
\usepackage{euscript}
\usepackage{epigraph}
\usepackage[numbers]{natbib}

\newcommand{\Z}{\mathbb{Z}}

\newcommand{\N}{\mathbb{N}}
\newcommand{\R}{\mathbb{R}}
\newcommand{\Q}{\mathbb{Q}}

\newcommand{\F}{\mathbb{F}}
\newcommand{\M}{\mathfrak{M}}

\newcommand{\g}{\mathfrak{g}}
\newcommand{\h}{\mathfrak{h}}

\newcommand{\eps}{\varepsilon}
\newcommand{\ph}{\varphi}

\newcommand{\ch}{\mathop{\mathrm{char}}\nolimits}

\newcommand{\sub}{\subseteq}

\newcommand{\Ker}{\mathop{\mathrm{Ker}}\nolimits}

\newcommand{\Spec}{\mathop{\mathrm{Spec}}\nolimits}

\newcommand{\rk}{\mathop{\mathrm{rk}}\nolimits}

\renewcommand{\ge}{\geqslant}
\renewcommand{\le}{\leqslant}
\newcommand{\sm}{\setminus}
\newcommand{\map}[3]{#1\colon #2\to #3}
\newcommand{\phan}{\vphantom{,}}

\newcommand{\restr}[2]{\left. #1\right|_{#2}}
\newcommand{\<}{\langle}
\renewcommand{\>}{\rangle}

\newcommand{\emp}{\varnothing}

\newcommand{\tc}{\text{,}}
\newcommand{\tp}{\text{.}}
\newcommand{\toiso}{\xrightarrow{\sim}}
\renewcommand{\tilde}{\widetilde}

\newcommand\leftact[2]{\phan^{#1}#2}

\DeclareMathOperator{\id}{id}

\DeclareMathOperator{\Pin}{Pin}
\DeclareMathOperator{\Adjust}{Adjust}
\DeclareMathOperator{\Isom}{Isom}
\DeclareMathOperator{\spl}{spl}
\DeclareMathOperator{\Lie}{Lie}
\DeclareMathOperator{\Diag}{Diag}
\DeclareMathOperator{\Aut}{Aut}
\DeclareMathOperator{\End}{End}
\DeclareMathOperator{\Der}{Der}
\DeclareMathOperator{\ad}{ad}

\theoremstyle{plain}
\newtheorem{thm}{Theorem}[section]
\newtheorem{lem}[thm]{Lemma}
\newtheorem{prop}[thm]{Proposition}
\newtheorem{cor}[thm]{Corollary}
\newtheorem*{vag}{Main Theorem stated vaguely}
\newtheorem*{exmp}{Example}

\theoremstyle{definition}
\newtheorem{defn}[thm]{Definition}

\theoremstyle{remark}
\newtheorem{rem}[thm]{Remark}

\AtBeginEnvironment{thm}{\begin{samepage}}
	\AtEndEnvironment{thm}{\end{samepage}}
\AtBeginEnvironment{lem}{\begin{samepage}}
	\AtEndEnvironment{lem}{\end{samepage}}
\AtBeginEnvironment{st}{\begin{samepage}}
	\AtEndEnvironment{st}{\end{samepage}}
\AtBeginEnvironment{crit}{\begin{samepage}}
	\AtEndEnvironment{crit}{\end{samepage}}
\AtBeginEnvironment{ax}{\begin{samepage}}
	\AtEndEnvironment{ax}{\end{samepage}}
\AtBeginEnvironment{defn}{\begin{samepage}}
	\AtEndEnvironment{defn}{\end{samepage}}
\AtBeginEnvironment{cor}{\begin{samepage}}
	\AtEndEnvironment{cor}{\end{samepage}}
\AtBeginEnvironment{note}{\begin{samepage}}
	\AtEndEnvironment{note}{\end{samepage}}
\AtBeginEnvironment{prop}{\begin{samepage}}
	\AtEndEnvironment{prop}{\end{samepage}}

\makeatletter
\def\@settitle{\begin{center}%
		\baselineskip14\p@\relax
		\bfseries
		\@title
	\end{center}%
}

\def\@evenhead{\hfil\sc pavel gvozdevsky\hfil}
\def\@oddhead{\hfil\sc abstract isomorphisms\hfil}
\makeatother

\title{Abstract isomorphisms of isotropic root graded groups over rings.}


\author{Pavel Gvozdevsky}
\date{}
\address{Department of Mathematics, Bar-Ilan University, 5290002 Ramat Gan, ISRAEL}
\thanks{The paper is written as part of the author's post-doctoral fellowship at Bar-Ilan University, Department of Mathematics; and is supported by the ISF grant 1994/20.
}
\email{gvozdevskiy96@gmail.com}

\begin{document}

\maketitle

\begin{abstract}
The celebrated Borel--Tits theorem provides a classification of abstract isomorphisms between (simple) isotropic groups over fields, showing that such isomorphisms arise from field isomorphisms and group-scheme isomorphisms. In this work, we extend the scope of this classification to certain class of group schemes over arbitrary commutative rings. Specifically, we prove that under suitable conditions abstract isomorphisms between the groups of points of isotropic, absolutely simple, adjoint group schemes over rings admit a description analogous to that in the classical setting: namely, they are induced by isomorphisms of ground rings and isomorphisms of the underlying group schemes. This result generalizes the classical theory to a far broader algebraic context and confirms that the rigidity phenomena observed over fields persist over rings.
\end{abstract}

\section{Introduction}

The celebrated Borel--Tits theorem provides a fundamental description of abstract isomorphisms between isotropic simple algebraic groups over fields. This landmark result was established in their seminal 1973 paper~\cite{BorelTitsHomo}. Almost immediately, the theorem became an indispensable tool in the theory of algebraic groups, underpinning a wide range of applications. Its influence was recognized, for example, in a dedicated session of the Séminaire N. Bourbaki (Talk 435), and in the subsequent survey by R. Steinberg~\cite{SteinbergBorelTits}, which distilled the core ideas and implications of the Borel--Tits theory.

Further developments culminated in a 2018 paper by A. Borel and J. Tits~\cite{BorelTitsHomoEng}, presenting a variation on the same theme. Despite these advances, one of the most significant and long-standing challenges in the area has been to extend the Borel--Tits framework beyond fields, particularly to arbitrary commutative rings and to anisotropic groups. This problem has remained open for decades and is widely regarded as both technically demanding and foundational for a deeper understanding of abstract homomorphisms in broader algebraic settings.

 More specifically, Theorem by Borel and Tits \cite[Theorem 8.11]{BorelTitsHomo} says that if $G_1$ and~$G_2$ are absolutely simple adjoint isotropic group schemes over infinite fields $K_1$ and~$K_2$, then any abstract isomorphism of their groups of points $G_1(K_1)\simeq G_2(K_2)$ must arise from a field isomorphism $\ph\colon K_1\toiso K_2$ and a $K_2$-group-scheme isomorphism $\Theta\colon \phan^\ph G_1\toiso G_2$ (with some exceptions). Actually, the Theorem says something a bit stronger, but for our purposes we only need this weaker version. The purpose of the present paper is to prove a similar statement for a broad class of group schemes over almost arbitrary commutative rings.

Such a result can be useful in many areas, in particular, it can play a decisive role in the study of the logic and model theory of isotropic reductive groups, cf.~\cite{BuninaGvoz}, ~\cite{KRT}.

To put the results of the present paper into historical context, we recall  the main landmarks  of research. For the general picture the surveys by Weisfeler and coauthors~\cite{WeisfeilerSurvey} and~\cite{HJW} are quite useful. The paper~\cite{JamesWaterhouseWeisfeiler} contains  a collection of problems. An updated bibliography can be found in the papers by I.Rapinchuk~\cite{Rapinchuk15} and Chatterjee~\cite{Chatterjee}. 

\smallskip

The first category of  results consists of those that are directly used in the present paper. 

$\bullet$ As mentioned above, the case of isotropic groups over infinite fields is covered in~\cite{BorelTitsHomo} (see also~\cite{SteinbergBorelTits} or \cite{BorelTitsHomoEng}); in fact, the paper describes not only isomorphisms, but also all the homomorphisms with dense image.

$\bullet$ The case of groups over finite fields is covered by a combination of the results in \cite{SteinbergAutFinite} and the knowledge of exceptional isomorphisms in the list of finite simple groups.

\smallskip

The second set of references consists of results that are substantially related to the ones of the present paper.

$\bullet$ Isomorphisms of split simple algebraic groups, that is of  Chevalley groups, were described by Steinberg (finite fields)~\cite{SteinbergAutFinite} and Humphreys~\cite{HumphreysAut}.

$\bullet$  Isomorphisms  of  Chevalley groups over rings are described in papers \cite{BuninaDifferentTypes,BuninaIsomorph,BuninaAllAut} by Bunina (see also \cite{AbeAut} by Abe, though the last paper contains a gap). These results for split groups essentially form a particular case of those in the present paper, except that Bunina assumes the existence of $1/2$ only for certain root systems, while we assume it in all cases. Also \cite{BuninaAllAut} describes a method to transfer results from adjoint groups to arbitrary ones, but it should be reconsidered  in the case of half-spinor groups. We also note that the methods of the present paper are largely inspired by those of Bunina.

$\bullet$ In papers \cite{PetechukMatrix,PetechukRk2,PetechukSp}, Petechuk describes automorphisms of general linear, special linear, and symplectic groups over a commutative ring. These results intersect with Bunina's; however, Petechuk allows rings without $1/2$, which are not covered by either Bunina’s papers or the present one. Apparently, in such cases the groups $\mathrm{SL}_3(R)$ and $\mathrm{GL}_3(R)$ may possess non-standard automorphisms.

$\bullet$ In \cite{LandinReiner} and \cite{ShijianYan}, automorphisms of the general linear group over not necessarily commutative principal ideal domains are described.

$\bullet$ In \cite{GolubchikIsoGen,GolMikhIsoUni,GolMikhIsoGen,Zelmanov}, isomorphisms of general linear groups and certain class of unitary groups over associative rings are treated.

\smallskip

The next results  develop the ideas of Borel and Tits in various directions.

$\bullet$ In \cite{Caprace}, isomorphisms of Kac--Moody groups over fields are described.

$\bullet$ In \cite{SeitzHomo}, arbitrary homomorphisms between Chevalley groups over fields are described.

$\bullet$ In \cite{WeisfeilerReal}, homomorphisms with dense image between anisotropic groups over real-closed fields are described. In~\cite{WeisfeilerQuadr}, similar results are obtained for arbitrary fields and for groups that split over a quadratic extension.

$\bullet$ A series of papers by I. Rapinchuk studies abstract homomorphisms between algebraic groups in certain cases where either the groups themselves or the image of the homomorphism may not be reductive. These results are under the flavor of Borel--Tits Conjecture (see~\cite[Remark 8.19]{BorelTitsHomo}). Its formulation is tightly related to the above discussed \cite[Theorem 8.11]{BorelTitsHomo} and despite the fact that the conjecture is wrong in general (see~\cite[Proposition 9.9.2]{CGP}), the series of papers \cite{Chatterjee,Rapinchuk11,LifschitzRapinchuk,Rapinchuk13,Rapinchuk19,RapinchukRuiter} sheds much light on abstract homomorphisms of algebraic groups. 

\smallskip

Finally, here are some references to earlier results that are special cases of those mentioned above:~\cite{HumphreysAut,SchreierWaerden,Hua51,Hua48,Diedonne,Rickart,OMeara66,OMeara68,McDonaldPomfret,WaterhouseGL,McQueenMcDonald,Klyachko,ChenYu95,ChenYu96,ChenYu00,CarterChen}.
Of course, this list is far from being complete.
\smallskip

The most common notion of an isotropic reductive group scheme over a ring is as follows: a reductive group scheme $G$ is said to be isotropic with isotropic rank $\ge n$, where $n \ge 1$, if $G$ contains a split torus of rank $n$ as a closed subgroup-subscheme. This is the definition used, for example, in \cite{PetStavIso} and \cite{StavStep}. However, we were not able to prove the theorem on abstract isomorphisms in this generality; instead, we require $G_1$ to be isotropic in the more restrictive sense of~\cite{VoronDiophantine}. This means that the following requirements, which are not automatic in general, must be satisfied: the relative root system in the sense of~\cite{PetStavIso} must actually be a root system; the corresponding root subgroups must make the groups $G_1(S)$ into a root graded group in the sense of~\cite{WiedemannRootGraded} and \cite{VoronRootGraded} for any $R$-algebra $S$ (alternatively, one can say that $G_1$ admits an isotropic pinning in the sense of~\cite{VoronIsotrStenberg}); and the resulting map from the absolute root system to the relative one must come from one of the Tits indexes. We believe it is natural to call this class of group schemes root graded isotropic group schemes and to refer to their groups of points as root graded isotropic groups.

The class of root graded isotropic group schemes includes, but is not limited to, special linear groups over Azumaya algebras; orthogonal groups corresponding to quadratic forms with positive Witt index; and unitary groups over Azumaya algebras with involution, associated with (anti)hermitian forms with positive Witt index. However, for example, if~$P$ is a projective module over the ring $R$, then a decomposition $P=P_1\oplus\ldots\oplus P_{n+1}$, where each $P_i$ has constant rank, leads to $\Aut_R(P)$ being a root graded isotropic group only if all the $P_i$ are isomorphic to each other, whereas $\Aut_R(P)$ is always isotropic in the usual sense, even if the ranks of $P_i$ are different from each other. Nevertheless, it follows from the results of \cite{PetStavTits} that over a semilocal ring all isotropic group schemes are root graded.

Additionally, for the sake of simplicity, we restrict ourselves to the case where the group schemes are not just reductive, but also absolutely simple and adjoint. The case of group schemes that are not necessarily adjoint will 
be treated separately in the future paper.

We now give the vague statement of our main theorem and one illustrative example; the precise statement will be given in Section~\ref{StatementSec}, all the terms involved will be explained in Section~\ref{BasicNotationSec}.

{\vag Let $G_1$ and $G_2$ be absolutely simple adjoint group schemes over rings $R_1$ and $R_2$, that satisfy certain technical assumptions {\rm(}including that~$G_1$ is root graded and that both have isotropic rank at least 2{\rm )}.

Let $E_i(R_i)$ be the elementary subgroup of $G_i(R_i)$ {\rm(}$i=1$,$2${\rm )}. Let $\theta\colon E_1(R_1)\toiso E_2(R_2)$ be the isomorphism of abstract groups.

Then there exists an isomorphism of rings $\ph\colon R_1\toiso R_2$ and an \mbox{$R_2$-group-scheme} isomorphism $\Theta\colon \phan^\ph G_1\toiso G_2$ such that $\theta=\restr{(\Theta_{R_2}\circ \ph_*)}{E_1(R_1)}$, where $\ph_*\colon  G_1(R_1)\toiso \phan^\ph G_1(R_2)$ is the isomorphism induced by $\ph$.
}

{\exmp Let $A_1$ and $A_2$ be two algebras such that each is an Azumaya algebra over its center. Suppose that 2 is invertible in the center of $A_2$. Let $n\ge 2$ and $m\ge 2$ be integers.

Suppose that elementary subgroups of the groups $\mathrm{PGL}_n(A_1)$ and $\mathrm{PGL}_m(A_2)$ are isomorphic. Then there are central idempotents $e_1\in A_1$ and $e_2\in A_2$ such that the matrix rings $M_{n\times n}(e_1A_1)$ and $M_{m\times m}(e_2A_2)$ are isomorphic; and the matrix rings $M_{n\times n}((1-e_1)A_1)$ and $M_{m\times m}((1-e_2)A_2)$ are anti-isomorphic.}

\smallskip

The paper is organized as follows. In Section~\ref{BasicNotationSec}, we introduce the basic notation and definitions, and explain everything required to understand the statement of the main theorem (Theorem~\ref{main}). In Section~\ref{StatementSec} we give the precise statement of the main theorem and sketch the main steps of the proof. Section~\ref{AuxilarySec} contains all the auxiliary stuff. In Sections~\ref{FieldSec}--\ref{ProofFinaleSec}, we present the proof of Theorem~\ref{main}. In Section~\ref{CorollarySec}, we show that the elementary subgroups in Theorem~\ref{main} can be replaced by any subgroups containing elementary ones. In Section~\ref{AssumptiondSec}, we explain the meaning of the assumption (d) in the second bullet in the statement of Theorem~\ref{main}.

The author is deeply grateful to Eugene Plotkin and Boris Kunyavskii for their constant support and interest in this work, and to Egor Voronetsky for sharing an early draft of~\cite{VoronDiophantine} before it appeared on the arXiv. I am  indebted to Nikolai Vavilov for introducing me the world of algebraic groups.

\section{Basic definitions, conventions and notation}
\label{BasicNotationSec}

\subsection*{Rings and algebras}$\phan$

The word ''ring'' always means a commutative associative ring with unity. If $R$ is a ring, the term ''$R$-algebra'' means a commutative associative $R$-algebra with unity, unless specified otherwise (e.g., when Azumaya algebras are mentioned).

If $R$ is a ring, we denote by $R^*$ the group of its invertible elements.

\subsection*{Operations in groups}$\phan$

We use the standard notation for conjugation in groups: $h^g=g^{-1}hg$ and $\leftact{g}{h}=ghg^{-1}$. Commutators are left-normalized: $[x,y]=xyx^{-1}y^{-1}$.

The notation $[\cdot,\cdot]$ may also refer to the Lie bracket in a Lie algebra, but the context will make it clear whether a calculation is performed in a group or an algebra.

\subsection*{Root systems}$\phan$

We usually use the notation $\tilde{\Phi}$ for a reduced irreducible crystallographic root system, and $\Phi$ for an irreducible crystallographic root system that is not necessarily reduced (i.e., it may be of type $BC$).

For roots $\alpha$ and $\beta$ we write $\<\beta,\alpha\>=\tfrac{2(\beta,\alpha)}{(\alpha,\alpha)}$, where $(\cdot,\cdot)$ denotes the inner product; and denote by $s_{\alpha}(\beta)$ the reflection of $\beta$ with respect to the hyperplane orthogonal to $\alpha$.

The root $\alpha$ is called {\it long} if $(\alpha,\alpha)\ge(\beta,\beta)$ for any root $\beta$ (in particular, in simply-laced systems every root counts as long). The root $\alpha$ is called {\it ultrashort} if $2\alpha$ is a root. The root~$\alpha$ is called {\it short} if it is neither long, nor ultrashort.

\subsection*{Schemes}$\phan$

We adopt the functorial point of view on the theory of schemes (see \cite{DemazureGabriel}). Thus, if $R$ is a ring, instead of saying that certain functor $X$ from the category of $R$-algebras to the category of sets is represented by a scheme over $R$, we will say that $X$ {\it is} a scheme over~$R$.

When we need to define a scheme over $R$, we specify what the set $X(S)$ is for each $R$-algebra $S$, and we promise to do so in a way that makes the implied action of $X$ on homomorphisms clear. The scheme in a more conventional sense (i.e., a locally ringed topological space) representing the functor $X$ will be referred to as the underlying topological space of the scheme $X$.

A morphism of schemes, in this terminology, is a natural transformation of functors. If $\map{\Theta}{X}{Y}$ is a morphism of schemes over $R$, and $S$ is an $R$-algebra we denote by $\Theta_S$ the corresponding map $X(S)\to Y(S)$.

A group scheme over $R$, in this terminology, is a functor from the category of $R$-algebras to the category of groups, whose composition with the forgetful functor to the category of sets is a scheme.

If $G$ is a group scheme over the ring $R$, and $S$ is an $R$-algebra with an ideal $I \unlhd S$, we denote by $G(S,I)$ the principal congruence subgroup of level $I$, i.e., the kernel of the reduction homomorphism $\map{\rho_I}{G(S)}{G(S/I)}$.

The notation $G(\tilde{\Phi},-)$ is reserved for the adjoint Chevalley--Demazure scheme associated with a reduced irreducible crystallographic root system $\tilde{\Phi}$. The Chevalley--Demazure schemes are defined in~\cite{ChevalleySemiSimp}. For a ring $S$ and an ideal $I \unlhd S$, we denote by $G(\tilde{\Phi},S)$ the corresponding adjoint Chevalley group, and by $G(\tilde{\Phi},S,I)$ the corresponding principal congruence subgroup of level~$I$.

\subsection*{Absolutely simple adjoint group schemes with a common root datum of the geometric fibers}$\phan$

For the purposes of the present paper, we do not need to define each term in the phrase “absolutely simple adjoint group scheme with a common root datum of the geometric fibers” separately; instead, we explain their combined meaning.

The group scheme $G$ over the ring $R$ is said to be an absolutely simple adjoint with a common root datum of the geometric fibers if there is an fppf-extension $S$ of the ring $R$ such that $G_S$ is isomorphic as a group scheme over $S$ to the adjoint Chevalley--Demazure scheme $G(\tilde{\Phi},-)_S$ over $S$ for some reduced irreducible crystallographic root system $\tilde{\Phi}$.

Here $G_S$ is a base change of $G$ to $S$. The root system $\tilde{\Phi}$ here is called the {\it absolute root system} of $G$.

\subsection*{Isotropic rank}$\phan$

We say that the absolutely simple adjoint group scheme $G$ with a common root datum of the geometric fibers has isotropic rank at least $l$ if it contains a closed subgroup-subscheme isomorphic to $\mathbb{G}_m^l$, where $\mathbb{G}_m$ is the multiplicative group scheme (i.e. $\mathbb{G}_m(S)=S^*$ for any $R$-algebra $S$).

\subsection*{Pinnings and isotropic pre-pinnings.}$\phan$

Let $G$ be an absolutely simple adjoint group scheme with a common root datum of the geometric fibers over the ring $R$. Recall \cite[Exp. XXIII, Definition 1.1]{SGAIII} that a \underline{pinning} of $G$ consists of

$\bullet$ a maximal torus $T\le G$ with a chosen isomorphism $T\simeq \mathbb{G}_m^l$;

$\bullet$ a root datum $(\Z^l,\Phi,(\Z^l)^\vee,\Phi^\vee)$ such that $\Phi$ and $\Phi^\vee$ are the sets of roots and coroots of $G$ with respect to $T$ (or $R=0$), in particular, roots and coroots are constant functions on $\Spec(R)$;

$\bullet$ a basis $\Delta\le\Phi$;

$\bullet$ trivializing sections $x_{\alpha}\in\g_{\alpha}$ in the root spaces $\g_{\alpha}\le\g=\Lie(G)$ for $\alpha\in\Delta$, so that all root spaces are free $R$-modules of rank 1.

\medskip

So if $G$ has a pinning, then it is isomorphic to $G(\tilde{\Phi},-)_R$ already over $R$ and the root system $\Phi$ in the corresponding root datum coincides with the absolute root system $\tilde{\Phi}$.

\medskip

Now following \cite{VoronIsotrStenberg}, we say that an \underline{isotropic pre-pinning} of $G$ consists of

$\bullet$ a split torus $T\le G$ with a chosen isomorphism $T\simeq \mathbb{G}_m^l$;

$\bullet$ a root system $\Phi \sub \Z^l$ (with respect to some inner product on $\R^l$) together with a
chosen base  $\Delta\le\Phi$ such that $\Phi$ is the set of non-zero weights of $\g$ with
respect to $T$ (or $R=0$);

\medskip

Here unlike \cite{VoronIsotrStenberg} we will not require the root subspaces to be free modules.

We say that an isotropic pre-pinning $(T, \Phi)$ is contained in an isotropic pre-pinning $(T',\Phi')$ if $T \le T'$, the inclusion is given by a constant surjective homomorphism $\map{u}{\Z^l}{\Z^l}$ of the corresponding abelian groups, and $u(\Phi'\cup\{0\})=\Phi\cup\{0\}$ (the last condition is vacuous if $R$ is non-zero).

For every isotropic pre-pinning $(T,\Phi)$ there is an fppf-extension $S$ of the ring $R$ such that after a base change to $S$ the isotropic pre-pinning $(T,\Phi)$ becomes contained in a pinning; therefore, inducing a map $\map{u}{\tilde{\Phi}\cup\{0\}}{\Phi\cup\{0\}}$. This map will be referred to as the corresponding map $\map{u}{\tilde{\Phi}\cup\{0\}}{\Phi\cup\{0\}}$ of the isotropic pre-pinning $(T,\Phi)$.

We will usually make an assumption that the corresponding map $\map{u}{\tilde{\Phi}\cup\{0\}}{\Phi\cup\{0\}}$ comes from one of the Tits indexes. The Tits index in question is always implied to be irreducible, because $\tilde{\Phi}$ is assumed to be irreducible. Recall that an irreducible Tits index $(\Phi, \Gamma, J)$ consists of a reduced irreducible crystallographic root system $\tilde{\Phi}$, a subgroup $\Gamma$ of its group of outer automorphisms (i.e. the automorphism group of its Dynkin diagram), and a $\Gamma$-invariant subset $J$ of vertices of the Dynkin diagram satisfying an additional condition (namely, that it may be constructed by a reductive group scheme over a field using its minimal parabolic subgroup). Then $\Phi$ is the image of $\tilde{\Phi}$ in the quotient-space of the ambient vector space by~$\Gamma$ and the span of basic roots not in $J$.

\subsection*{Root subgroups and elementary subgroup}$\phan$

For every isotropic pre-pinning there exist root subgroup-subschemes $G_{\alpha}\le G$ uniquely defined by certain list of properties. If $\alpha\in\Phi$ is not ultrashort, then $G_{\alpha}(S)$ is an abelian group with a natural structure of $S$-module. In case where $\alpha$ is ultrashort, we have $G_{2\alpha}(S)\le G_{\alpha}(S)$ and there is a natural action $(-)\cdot (=)$ of the multiplicative monoid $S^\bullet$ of the ring $S$ on $G_{\alpha}(S)$ by group endomorphisms from the right, which induces an $S$-module structure on the abelian quotient group $G_{\alpha}(S)/G_{2\alpha}(S)$. For the details see \cite{VoronLocIsotr} or \cite{VoronIsotrStenberg}.

The subgroup $E(S)\le G(S)$ generated by all the $G_{\alpha}(S)$ is called the {\it elementary subgroup}. More generally, the elementary subgroup is defined in~\cite{PetStavIso} for isotropic reductive groups as a subgroup generated by subgroups of points of two opposite unipotent radicals.

\subsection*{Weyl elements and isotropic pinnings}$\phan$

Let $G$ be an absolutely simple adjoint group scheme with a common root datum of the geometric fibers over the ring $R$; and let $(T,\Phi)$ be an isotropic pre-pinning on $G$. We adopt the following definitions from \cite{VoronIsotrStenberg}.

For a root $\alpha\in\Phi$, an element $w_{\alpha}\in G_{\alpha}(R) G_{-\alpha}(R) G_{\alpha}(R)$ is called a {\it Weyl element} if $\leftact{w_{\alpha}}{G_{\beta}(S)}=G_{s_{\alpha}(\beta)}(S)$ for every root $\beta\in\Phi$ and any $R$-algebra $S$.

An isotropic pre-pinning $(T,\Phi)$ is called an {\it isotropic pinning} if Weyl elements $w_{\alpha}$ exist for every $\alpha\in\Phi$.

\section{Statement of the main result and sketch of the proof}
\label{StatementSec}

We now give the precise statement of our main theorem.

{\thm\label{main} Let $G_1$ and $G_2$ be absolutely simple adjoint group schemes over rings $R_1$ and $R_2$, each with a common root datum of the geometric fibers. Let $\tilde{\Phi_1}$ and $\tilde{\Phi_2}$ be the corresponding absolute root systems. Assume the following:
	
	$\bullet$ $G_2$ has isotropic rank at least 2;
	
	$\bullet$ $G_1$ admits an isotropic pinning with the following properties:
	
	a) its root system $\Phi$ has rank at least 2;
	
	b) the corresponding map $\map{u}{\tilde{\Phi_1}\cup\{0\}}{\Phi\cup\{0\}}$ comes from one of the Tits indexes;
	
	c) it has square formula {\rm(}see Definition~\ref{SquareFormula}{\rm)};
	
	d)if $\Phi$ is of type $C$ or $BC$ (including $C_2$) and the map $u$ is not a bijection, then for every pair of orthogonal short roots $\beta$,$\beta'\in\Phi$ with their sum being a long root the corresponding Weyl elements $w_{\beta}=a_{\beta}b_{\beta}c_{\beta}$ and $w_{\beta'}=a_{\beta'}b_{\beta'}c_{\beta'}$ {\rm(}here $a_{\beta}$ and $c_{\beta}$ belong to the root subgroup of the root $\beta$, while $b_{\beta}$ belong to the root subgroup of the root $-\beta$, and similarly for $a_{\beta'}$, $b_{\beta'}$ and $c_{\beta'}${\rm)} can be chosen so that $b_{\beta}$ commutes with $b_{\beta'}$;
	
	$\bullet$ if $\tilde{\Phi_1}$ is doubly laced, then $2\in R_1^*$; if $\tilde{\Phi}_1=G_2$, then $6\in R_1^*$;
	
	$\bullet$ $2\in R_2^*$ (regardless of $\tilde{\Phi_2}$); if $\tilde{\Phi_2}=G_2$, then $6\in R_2^*$.
	
	Let $E_i(R_i)$ be the elementary subgroup of $G_i(R_i)$ {\rm(}$i=1$,$2${\rm )}. Let $\theta\colon E_1(R_1)\toiso E_2(R_2)$ be the isomorphism of abstract groups.
	
	Then
	
	\begin{enumerate}
		\item If $\tilde{\Phi_1}$ is not isomorphic to $\tilde{\Phi_2}$, then $\tilde{\Phi_1}=A_3$, $\tilde{\Phi_2}=B_2$, $R_1/\M\simeq \F_2$ for all maximal ideals $\M\unlhd R_1$ and $R_2/\M\simeq \F_3$ for all maximal ideals $\M\unlhd R_2$.
		
		\item  If $\tilde{\Phi_1}=\tilde{\Phi_2}$, then there exists an isomorphism of rings $\ph\colon R_1\toiso R_2$ and an \mbox{$R_2$-group-scheme} isomorphism $\Theta\colon \phan^\ph G_1\toiso G_2$ such that $\theta=\restr{(\Theta_{R_2}\circ \ph_*)}{E_1(R_1)}$, where $\ph_*\colon  G_1(R_1)\toiso \phan^\ph G_1(R_2)$ is the isomorphism induced by $\ph$.
	\end{enumerate}
}

We now give a brief overview of the proof.

In Section~\ref{MaximalSec} we use the results of \cite{StavStep} in order to define a one-to-one correspondence between the maximal ideals of $R_1$ and $R_2$ Such that $\theta$ induces an isomorphism $E_1(R_1/I)\simeq E_2(R_2/J)$ whenever $I$ corresponds to $J$. Note that $R_1/I$ and $R_2/J$ are fields.

The isomorphisms of elementary groups over fields can be described with some exceptions using results of \cite{BorelTitsHomo},\cite{SteinbergAutFinite} and the knowledge of coincidences in the list of finite simple groups, see Section~\ref{FieldSec} for details. This immediately imply the item (1) of Theorem~\ref{main}, see Section~\ref{MaximalSec}.

One of the steps we do for the proof of item (2) is to construct an isotropic pinning on the scheme $G_2$  with the same Tits index as the isotropic pinning on $G_1$ such that the isomorphism $\theta$ maps root subgroups to the corresponding root subgroups. This step is done in Section~\ref{PinningSec} using a descent from a faithfully flat extension $S$ of the ring $R_2$. In order for this to work the extension $S$ must be constructed in such a way that certain scheme, which we call the scheme of adjustments (it essentially describes pinnings on $G_2$ with some additional conditions needed for the analysis of the isomorphism $\theta$, see Section~\ref{AdjustmentSec}) has a point over $S$.

In order to construct such extension, in Section~\ref{AdjustmentSec} we will prove that the structural morphism from the scheme of adjustment to $\Spec(R_2)$ is smooth and surjective on the underlying topological spaces. Once we prove smoothness, the surjectivity will easily follow from the results of Sections~\ref{MaximalSec} and~\ref{FieldSec}; however, proving smoothness is not easy. We will use the following criterion from \cite[Ch I, \S 4, Item 4.5]{DemazureGabriel}: suppose that the scheme $X$ is locally finitely presented over the ring $R$, then $X$ is smooth over $R$ iff for ny $R$ algebra $S$ and any ideal $I\unlhd S$ with $I^2=0$ the reduction map $X(S)\to X(S/I)$ is surjective. In order to apply this criterion to the scheme of adjustments we will need a technical result from Section~\ref{InfinitesimalSecI} involving lengthy computations with root subgroups and commutator relations; so the umbrella assumptions in this section are desined to address this specific problem.

Once we construct the suitable isotropic pinning on $G_2$, we use the results of~\cite{VoronDiophantine} in order to construct a ring isomorphism $\ph$ between $R_1$ and $R_2$, see Section~\ref{RingsSec}. After that, in Section~\ref{IsomSec}, we define the scheme $\Isom^\#$ over the ring $R_2$, which describes group-scheme isomorphisms between $\phan^\ph G_1$ and $G_2$ that coincide with $\theta$ on the elementary subgroup. So in order to finish the proof of item (2) of Theorem~\ref{main}, we must show that the scheme $\Isom^\#$ has a point over $R_2$.

The three main steps of the remaining proof are: proving that the scheme $\Isom^\#$ has at most one point over any $R_2$ algebra, proving that the scheme $\Isom^\#$ is smooth over~$R_2$, and proving that its structural morphism to $\Spec(R_2)$ is surjective on the underlying topological spaces. The main difficulty here is the proof of smoothness, which will require us to reuse the results of Section~\ref{InfinitesimalSecI}, and also to use certain result from Section~\ref{AuxilarySec} on ''diagonal'' automorphisms.

\section{Auxiliary stuff}
\label{AuxilarySec}

\subsection*{Weyl elements in the split case}

{\lem\label{WeylForSplit} Let $\map{u}{\tilde{\Phi}\cup\{0\}}{\Phi\cup\{0\}}$ be the map of root systems that comes from one of the Tits indexes. Consider an isotropic pre-pinning $(T,\Phi)$ on the adjoint Chevalley--Demazure scheme $G(\tilde{\Phi},-)$ that is contained in the pinning and such that the corresponding map ${\tilde{\Phi}\sqcup\{0\}}\to{\Phi\sqcup\{0\}}$ is $u$. Then $(T,\Phi)$ is actually an isotropic pinning, i.e. there are Weyl elements with respect to $\Phi$.
}
\begin{proof}
	Clearly, it is enough to prove the existence of the Weyl elements $w_{\alpha}$ for all non-ultrashort roots $\alpha\in\Phi$. Now let $\alpha\in\Phi$ be a non-ultrashort root. The case-by-case consideration of Tits indexes shows that there exists a non-empty set of roots $\{\tilde{\alpha}_i\}\sub u^{-1}(\alpha)$ such that all the $\tilde{\alpha}_i$ are of the same length, pairwise orthogonal (and the sum of any pair is not a root), and the sum of all the $\tilde{\alpha}_i$ is orthogonal to the kernel of $u$ considered as a map of the ambient vector spaces. Now it is easy to see that we can take $w_\alpha=\prod w_{\tilde{\alpha}_i}$.
	\end{proof}

\subsection*{Elements $h_{\alpha}$ and the square formula}

{\lem\label{haExist} Let $G$ be an absolutely simple adjoint group scheme over a ring $R$ with an isotropic pre-pinning $(T,\Phi)$. Suppose that the corresponding map $\map{u}{\tilde{\Phi}\cup\{0\}}{\Phi\cup\{0\}}$ comes from one of the Tits indexes. Then for any non-ultrashort $\alpha\in\Phi$ there exists a unique element $h_{\alpha}\in G(R)$ such that for every root $\beta\in\Phi$, every $R$-algebra $A$, and any element $u\in G_{\beta}(A)$ we have
$$
u^{h_{\alpha}}=\begin{cases}
	u\;&\text{ if } \<\beta,\alpha\>\text{ is even }\\
	u^{-1}\;&\text{ if } \<\beta,\alpha\>\text{ is odd and }\beta\text{ is non-ultrashort}\\
	u\cdot(-1)\;&\text{ if } \<\beta,\alpha\>\text{ is odd and }\beta\text{ is ultrashort}\tp
	\end{cases}
$$}
\begin{proof}
Let $S$ be a faithfully flat extension of $R$ such that the isotropic pinning $(T,\Phi)$ on $G$ after a base change to $S$ becomes contained in a pinning. It is easy to see that such $h_{\alpha}$ exist over $S$. These $h_{\alpha}$ are uniquely determined by the required property, and this definition uses only the subschemes $G_{\beta}$, which are defined over $R$. Thus the maps $\id\otimes 1$ and $1\otimes \id$ from $S$ to $S\otimes_R S$ must coincide on the elements $h_{\alpha}\in G(S)$. Therefore, by flat descent we obtain that $h_{\alpha}$ must exist over $R$.
	\end{proof}

{\defn\label{SquareFormula} We say that an isotropic pinning $(T,\Phi)$ has {\it square formula}, if the Weyl elements $w_{\alpha}$ can be chosen so that $w_{\alpha}^2=h_{\alpha}$ for all $\alpha\in\Phi$.}

{\rem The square formula is known to be true for the cases, where either $\Phi$ is simply-laced or $\rk\Phi\ge 3$ see \cite[Propositions 5.4.15; 7.6.15 and 9.5.13]{WiedemannRootGraded}. Hypothesis: it is always the case.}

\subsection*{Defining root subgroups by formulas}

{\lem\label{EgorFormulas} Let $G$ be an absolutely simple adjoint group scheme with a common root datum of the geometric fibers over ring $R$ with an absolute root system $\tilde{\Phi}$. Let $(T,\Phi)$ be an isotropic pre-pinning on $G$ such that $\rk\Phi\ge 2$ and such that the corresponding map $\map{u}{\tilde{\Phi}\cup\{0\}}{\Phi\cup\{0\}}$ comes from one of the Tits indexes. Let $G_{\alpha}$ for $\alpha\in\Phi$ be the corresponding root subgroup-subschemes. For each $\alpha\in\Phi$ let $X_{\alpha}\sub G_{\alpha}(R)$ be a subset that generate $G_{\alpha}(R)$ (or $G_{\alpha}(R)/G_{2\alpha}(R)$ if $\alpha$ is ultrashort) as an $R$-module. Then for every $R$-algebra $S$ and for any $\alpha\in\Phi$ we have
$$
G_{\alpha}(S)=\{g\in G(S)\colon \text{(F1)--(F3) hold}\}\tc
$$
where

(F1) $[g,x_{\gamma}]=e$ for every $\gamma\in\Phi$ such that $\alpha+\gamma\notin\Phi\cup\{0\}$ and every $x_{\gamma}\in X_{\gamma}$.

(F2) If $\alpha$ is ultrashort, then $[[g,x_{\alpha}],x_{\gamma}]=e$ for every $\gamma\in\Phi$ such that $2\alpha+\gamma\notin\Phi\cup\{0\}$ and every $x_{\gamma}\in X_{\gamma}$ and $x_{\alpha}\in X_{\alpha}$.

(F3) If $\alpha$ is a short root and $\Phi$ is of type $C$ or $BC$ (including $C_2=B_2$), and $\alpha'$ and~$\alpha''$ are distinct long roots such that $\angle(\alpha,\alpha')=\angle(\alpha,\alpha'')=\pi/4$, then $[[g,x_{\alpha'-\alpha}],x_{\gamma}]=e$ for every $\gamma\in\Phi$ such that $\alpha'+\gamma\notin\Phi\cup\{0\}$ and every $x_{\gamma}\in X_{\gamma}$ and $x_{\alpha'-\alpha}\in X_{\alpha'-\alpha}$; and $[[g,x_{\alpha''-\alpha}],x_{\gamma}]=e$ for every $\gamma\in\Phi$ such that $\alpha''+\gamma\notin\Phi\cup\{0\}$ and every $x_{\gamma}\in X_{\gamma}$ and $x_{\alpha''-\alpha}\in X_{\alpha''-\alpha}$.
}
\begin{proof}
	Since the statement is fppf-local, by Lemma~\ref{WeylForSplit} we may assume that $(T,\Phi)$ is an isotropic pinning.
	
	Now if $\alpha$ is neither ultrashort, nor a short root in a system of type $C$ or $BC$, then the statement follows directly from \cite[Theorem 2, item 3]{VoronDiophantine}.
	
	If $\alpha$ is ultrashort, then by \cite[Theorem 2, item 1]{VoronDiophantine} F1 is equivalent to $g\in C^{\mathrm{us}}G_{\alpha}(S)$, where $C^{\mathrm{us}}$ is the scheme centralizer of non-ultrashort root subgroups; and by the previous case F2 is equivalent to $[g,X_{\alpha}]\le G_{2\alpha}(S)$. Hence together they are equivalent to $g\in G_{\alpha}(S)$ (see the proof of \cite[Theorem 4]{VoronDiophantine}).
	
	If $\alpha$ is a short root and $\Phi$ is of type $C$ or $BC$, and $\alpha'$ and~$\alpha''$ are distinct long roots such that $\angle(\alpha,\alpha')=\angle(\alpha,\alpha'')=\pi/4$, then by \cite[Theorem 2, item 2]{VoronDiophantine} F1 is equivalent to $g\in G_{\alpha}(S)G_{\alpha'}(S)G_{\alpha''}(S)$; and by the first case F3 is equivalent to $[g.X_{\alpha'-\alpha}]\le X_{\alpha'}$ and   $[g.X_{\alpha''-\alpha}]\le X_{\alpha''}$. It now follows from \cite[Lemma 3]{VoronDiophantine} that together they are equivalent to $g\in G_{\alpha}(S)$.
\end{proof}

\subsection*{Automorphisms of an affine line as an algebraic ring}

{\lem\label{RingStays} Let $R$ be a ring. Then the forgetful functor from the category of $R$-algebras to the category of rings has no non-trivial automorphisms.}

\begin{proof}
Let $\eta$ be an automorphism of the forgetful functor. Set $f=\eta_{R[x]}(x)$ and $g=\eta_{R[x]}^{-1}(x)$. So for any $R$-algebra $A$ and any $a\in A$ we have $\eta_A(a)=f(a)$.

By construction we have
$$
f(g(x))=g(f(x))=x\eqno{(1)}\tp
$$
Since all the $\eta_A$ must be automorphisms of rings we have
$$
	f(x+y)=f(x)+f(y)\eqno{(2)}
	$$
	and
	$$
	f(xy)=f(x)f(y)\tp\eqno{(3)}
$$
Similarly, we have
$$
g(x+y)=g(x)+g(y)\tp\eqno{(2')}
$$

Comparing free coefficients in $(2)$ and $(2')$ we obtain that the free coefficients of $f$ and $g$ vanish. Taking this into account, comparing the coefficient in $x$ in $(1)$ we obtain that $f$ has invertible coefficient in $x$. Finally, comparing all the coefficients in $(3)$ we obtain that the coefficients of $f$ are pairwise orthogonal idempotents. Therefore, we have $f(x)=x$; hence, $\eta$ is trivial.
\end{proof}

\subsection*{Clopen embeddings of group schemes}$\phan$

\smallskip

The following Lemma is a variation of \cite[Theorem 1.6.1]{WaterhouseAut}

{\lem\label{Waterhouse} Let $\mathcal{S}$ be set of prime numbers. Let $G$ and H be affine group schemes of finite type over $\Z[\mathcal{S}^{-1}]$ with $G$ flat, and let $\map{\ph}{G}{H}$ be a group-scheme homomorphism. Suppose that for all algebraically closed fields $K$ with $\ch K\notin\mathcal{S}$ we have
\begin{enumerate}
\item $\dim (G_K)\ge \dim_K\Lie(H_K)$;

\item the maps $G(K)\to H(K)$ and $G(K[\eps]/(\eps^2))\to H(K[\eps]/(\eps^2))$ given by $\ph$ are injective.
\end{enumerate}
Then $\ph$ is simultaneously a closed embedding and an open embedding.
 }
\begin{proof}
	Let us prove that for any algebraically closed field $K$ with $\ch L\notin\mathcal{S}$ the homomorphism $\map{\ph_K}{G_K}{H_K}$ is a clopen embedding.
	
	First (2) implies that $\dim(G_K)\le \dim (H_K)$. Together with (1), this implies that $G$, and~$H$, are both smooth and have the same dimension. Now by \cite[Ch II, \S 5, Item 5.5 (a)]{DemazureGabriel} $\Ker \ph_K$ is etale. Together with the fact that by (2) $(\Ker \ph_K)(K)$ is trivial this implies that $\Ker \ph_K$ is trivial; hence, $\ph_K$ is a monomorphism. Therefore, $\ph_K$ is a closed embedding by  \cite[Ch II, \S 5, Item 5.1 (b)]{DemazureGabriel}; and an open embedding by \cite[Ch II, \S 5, Item 5.5 (b)]{DemazureGabriel}.
	
	Further, if $k$ is any field with $\ch k\notin\mathcal{S}$, then the descent from $K=\overline{k}$ to $k$ shows that $\ph_k$ is a clopen embedding.
	
	Now consider the map $\map{\ph^*}{\Z[\mathcal{S}^{-1}][H]}{\Z[\mathcal{S}^{-1}][G]}$ on the rings of regular functions on the schemes $H$ and $G$ induced by $\ph$. Let $M\unlhd \Z[\mathcal{S}^{-1}][H]$ be the kernel of $\ph^*$. The fact that $\ph_{\Q}$ is a clopen embedding implies that $\map{\ph_{\Q}^*}{\Z[\mathcal{S}^{-1}][H]\otimes\Q}{\Z[\mathcal{S}^{-1}][G]\otimes\Q}$ is a surjection with a kernel generated by an idempotent, which implies that $M/M^2$ is a torsion group. As $\Z[\mathcal{S}^{-1}][H]$ is a finitely generated $\Z[\mathcal{S}^{-1}]$-algebra, the ideal $M$ is a finitely	generated $\Z[\mathcal{S}^{-1}][H]$-module, and hence there is a bound on the orders of elements in $M/M^2$.
	
	We now prove that $M=M^2$. Suppose that for some prime $p\notin\mathcal{S}$ there is a nontrivial $p$-primary component of $M/M^2$. Then we can choose an element of largest $p$-power order, and it will give a nontrivial class in $M/(M^2+pM)$. The fact that $\ph_{\F_p}^*$ is a clopen embedding implies that $x=x'+px''$ for some $x'\in M$ and $x''\in \Z[\mathcal{S}^{-1}][H]$. This implies that $px''\in M$. By assumption $\Z[\mathcal{S}^{-1}][G]$ is a flat $\Z[\mathcal{S}^{-1}]$-algebra; hence it is torsion free. Therefore, $\ph^*(px'')=0$ implies $\ph^*(x'')=0$; i.e. we have $x''\in M$; hence, $x\in M^2+pM$. This is a contradiction that shows that $M=M^2$.
	
	Since $\Z[\mathcal{S}^{-1}][H]$ is a finitely generated $\Z[\mathcal{S}^{-1}]$-algebra, it follows by Nakayama's lemma that $M$ is generated by some idempotent $e\in\Z[\mathcal{S}^{-1}][H]$.
	
	We now show that $\ph^*$ is surjective. Suppose we could find some $x\in \Z[\mathcal{S}^{-1}][G]$ not in the image. We know that $\ph^*$ becomes a surjection over $\Q$, so by appropriate multiplication we can assume that for some prime $p\notin\mathcal{S}$ we have $px=\ph^*(y)=\ph^*((1-e)y)$. But $px$ becomes zero when we tensor with $\F_p$. Hence, we have $(1-e)y=ey'+py''$ for some $y'$,$y''\in \Z[\mathcal{S}^{-1}][H]$. Multiplying by $(1-e)$ we get $(1-e)y=p(1-e)y''$. Hence we have $p(x-\ph^*((1-e)y''))=px-\ph^*((1-e)y)=0$; and since $\Z[\mathcal{S}^{-1}][G]$ is torsion-free, it follows that $x=\ph^*((1-e)y'')$. This is a contradiction that shows that $\ph^*$ is surjective.
	
	Since  $\ph^*$ is surjective with a kernel generated by an idempotent, it follows that $\ph$ is a clopen embedding.
\end{proof}

\subsection*{Presentation of classical Lie algebras}$\phan$

\smallskip

Let $\tilde{\Phi}$ be a reduced irreducible crystallographic root system. Let $\g({\Z})$ be an integer span of the Chevalley basic in the simple Lie algebra of type $\tilde{\Phi}$; and for any ring $R$ set $\g(R)=g(\Z)\otimes_{\Z} R$. Alternatively, $\g(R)$ can be viewed as a tangent Lie algebra to the simply connected Chevalley--Demazure scheme over $R$.

Let $\map{u}{\tilde{\Phi}\cup\{0\}}{\Phi\cup\{0\}}$ be a map that comes from one of the Tits indexes. Let $\g_{\alpha}(R)$ for $\alpha\in\Phi$ be the sum of root subspaces for all roots from $u^{-1}(\alpha)$.

{\lem\label{LieCopresent} Suppose that $\rk\Phi\ge 2$. Let $k$ be a field with $\ch k\ne 2$. If $\tilde{\Phi}=G_2$ assume also that $\ch k\ne 3$. Then the algebra $\g=\g(k)$ can be defined as a Lie algebra over $k$ by generators and relations in the following way.

Generators: elements from $\g_{\alpha}$ for all $\alpha\in\Phi$.

Relations:

(R1) all the linear relations from spaces $\g_{\alpha}$;

(R2) if $\alpha+\beta\notin\Phi\cup\{0\}$, then $[x,y]=0$ for all $x\in\g_{\alpha}$, $y\in\g_{\beta}$;

(R3) if $\alpha$,$\beta\in\Phi$ with $\alpha+\beta\in\Phi$ are either independent or equal, then $[x,y]=z$ for all $x\in\g_{\alpha}$, $y\in\g_{\beta}$, where $z\in\g_{\alpha+\beta}$ is the actual Lie bracket of $x$ and $y$ in $\g$.}

\begin{proof}
	Let $\g^\dag$ be the Lie algebra defined by generators and relations as above. We will write $[x,y]$ for the Lie bracket in $\g$ and $[x,y]_{\dag}$ for the Lie bracket in $\g^\dag$. Let $\map{\tau}{\g^\dag}{\g}$ be the homomorphism that is identity on each $\g_{\alpha}$.
	
	\smallskip
	
	\underline{Sublemma A} For any $\gamma\in\Phi$ there are independent roots $\alpha$,$\beta\in\Phi$ with $\alpha+\beta=\gamma$ such that $[\g_{\alpha},\g_{\beta}]=\g_{\gamma}$ (which by R3 implies that $[\g_{\alpha},\g_{\beta}]_{\dag}=\g_{\gamma}$).
	
	\medskip
	
	\underline{Proof of Sublemma A.} In most of the cases we can find independent roots $\alpha$,$\beta\in\Phi$ with $\alpha+\beta=\gamma$ such that $\beta$ is not ultrashort and $\gamma=s_{\beta}(\alpha)$; hence, it follows from Lemma~\ref{WeylForSplit} that $[\g_{\alpha},\g_{\beta}]=\g_{\gamma}$.
	
	The only case, where we can not do as above, is the case where $\gamma$ is a long root and $\Phi$ has type $C$ or $BC$. In this case we can take $\alpha$,$\beta\in\Phi$ to be the short roots with $\alpha+\beta=\gamma$. It is easy to verify case by case, that we have $2\g_{\gamma}\sub [\g_{\alpha},\g_{\beta}]$; and since we have $\ch k\ne 2$, it follows that $[\g_{\alpha},\g_{\beta}]=\g_{\gamma}$.
	
	 \underline{Sublemma A is proved.}
	
	\medskip
	
	\underline{Sublemma B} If $\alpha$,$\beta\in\Phi$ are such that $\alpha+\beta\in\Phi$, then $[x,y]_{\dag}=[x,y]$ for all $x\in\g_{\alpha}$, $y\in\g_{\beta}$.
	
	\medskip
	
	\underline{Proof of Sublemma B.} If $\alpha$ and $\beta$ are either independent or equal, then this is just R3. Thus, we only have to consider the case, where $\Phi$ has type $BC$ and $\beta=-2\alpha$.
	
	Let $\alpha'$ be an ultrashort root and $\alpha''$ be a short root such that $\alpha'+\alpha''=\alpha$. It follows from Lemma~\ref{WeylForSplit} that $x=[x',x'']=[x',x'']_{\dag}$ for some $x'\in\g_{\alpha'}$ and $x''\in\g_{\alpha''}$. Hence, we have
	\begin{multline*}
	[x,y]_{\dag}=[[x',x'']_{\dag},y]_{\dag}=[[x',y]_{\dag},x'']_{\dag}+[x',[x'',y]_{\dag}]_{\dag}=[x',[x'',y]_{\dag}]_{\dag}=[x',[x'',y]]_{\dag}=\\=[x',[x'',y]]=[[x',x''],y]+[x'',[x',y]]=[x,y]\tp
\end{multline*}
	
	\underline{Sublemma B is proved.}
	
	\medskip
	
	Now let $\h^{\dag}\le \g^{\dag}$ be the intersection of normalizers of $\g_{\alpha}$ in $\g$ for all $\alpha\in\Phi$.
	
	\medskip
	
	\underline{Sublemma C} For any $\alpha\in\Phi$ we have $[\g_{\alpha},\g_{-\alpha}]_{\dag}\le\h^{\dag}$.
	
	\medskip
	
	\underline{Proof of Sublemma C.} We must prove that $[[\g_{\alpha},\g_{-\alpha}]_{\dag},\g_{\beta}]_{\dag}\le\g_{\beta}$ for all $\beta\in\Phi$. First, let us consider the case, where $\alpha$ and $\beta$ are independent. We have
	$$
	[[\g_{\alpha},\g_{-\alpha}]_{\dag},\g_{\beta}]_{\dag}=[[\g_{\alpha},\g_{\beta}]_{\dag},\g_{-\alpha}]_{\dag}+[\g_{\alpha},[\g_{-\alpha},\g_{\beta}]_{\dag}]_{\dag}\tp
	$$
	Let us prove that $[[\g_{\alpha},\g_{\beta}]_{\dag},\g_{-\alpha}]_{\dag}\le \g_{\beta}$. If $\alpha+\beta\notin\Phi$, then we have $[[\g_{\alpha},\g_{\beta}]_{\dag},\g_{-\alpha}]_{\dag}=0$ by R2. If $\alpha+\beta\notin\Phi$, then by R3, we have
	$$
	[[\g_{\alpha},\g_{\beta}]_{\dag},\g_{-\alpha}]_{\dag}\le[\g_{\alpha+\beta},\g_{-\alpha}]_{\dag}\le \g_{\beta}\tp
	$$
	Similarly, we have $[\g_{\alpha},[\g_{-\alpha},\g_{\beta}]_{\dag}]_{\dag}\le \g_{\beta}$.
	
	Now, if $\beta\in\R\alpha$, then by Sublemma A, we have $\g_{\beta}=[\g_{\beta'},\g_{\beta''}]_{\dag}$, where $\beta=\beta+\beta'$, and $\beta'$ and $\beta''$ are independent; hence each of them is independent with $\beta$; and hence with~$\alpha$. Therefore, by what we prove previously, $[\g_{\alpha},\g_{-\alpha}]_{\dag}$ normalizes $\g_{\beta'}$ and $\g_{\beta''}$; hence it normalizes $\g_{\beta}$.
	
	\underline{Sublemma C is proved.}
	
	\medskip
	
		\underline{Sublemma D} The algebra $\g^{\dag}$ is the sum of subspaces $\h^{\dag}$ and all the $\g_{\alpha}$.
	
	\medskip
	
	\underline{Proof of Sublemma D.} By definition $\g^{\dag}$ is generated by all the $\g_{\alpha}$ as an algebra. Hence it is enough to show that the sum of subspaces $\h^{\dag}$ and all the $\g_{\alpha}$ is a subalgebra. It is easy to see that $[h^{\dag},h^{\dag}]_{\dag}\le h^{\dag}$. By definition we have $[h^{\dag},g_{\alpha}]_{\dag}\le g_{\alpha}$. If $\alpha+\beta\notin\Phi\cup\{0\}$, then by R2 we have $[g_{\alpha},g_{\beta}]_{\dag}=0$. If $\alpha+\beta\in\Phi$, then by Sublemma B we have $[g_{\alpha},g_{\beta}]_{\dag}=\g_{\alpha+\beta}$. Finally, by Sublemma~C we have $[\g_{\alpha},\g_{-\alpha}]_{\dag}\le\h^{\dag}$.
	
	\underline{Sublemma D is proved.}
	
	\medskip
	
	\underline{Sublemma E} $\Ker\tau$ is contained in the center of $\g^\dag$.
	
	\medskip
	
	\underline{Proof of Sublemma E.} Let $\h\le\g$ be the sum of the Cartan subalgebra with all the root subspaces for roots from $u^{-1}(0)$. It is easy to see that $\tau(\h^{\dag})\le \h$; and by definition we have $\tau(\g_{\alpha})=\g_{\alpha}$.
	
	Now it follows from Sublemma D that $\Ker\tau\le\h^{\dag}$. Hence, we have $[\Ker\tau,\g_{\alpha}]\le\g_{\alpha}$ for all $\alpha\in\Phi$. Hence, we have
	$$
	[\Ker\tau,\g_{\alpha}]_{\dag}=\tau([\Ker\tau,\g_{\alpha}]_{\dag})=[\tau(\Ker\tau),\g_{\alpha}]=0\tp
	$$
	
	 Therefore, $\Ker\tau$ is contained in the center of $\g^\dag$.
	
	\underline{Sublemma E is proved.}
	
	\medskip
	
	Now we finish the proof of the lemma. Since $\ch k\ne 2$, it follows that $\g=[\g,\g]$. Hence there is a universal central extension $\map{\pi}{g^*}{g}$. It is easy to see that $\tau$ is surjective; hence, by Sublemma E we have $\pi=\tau\circ\pi^\dag$ for a unique homomorphism $\map{\pi^\dag}{g^*}{g^\dag}$. It follows from Sublemma A that $\g^{\dag}=[\g^{\dag},\g^{\dag}]_{\dag}$; hence, $\pi^\dag$ is surjective. Now it only remains to prove that $\Ker \pi\le\Ker\pi^\dag$.
	
	Since we have $\ch k\ne 2$ and, in case $\tilde{\Phi}=G_2$, we also have $\ch k\ne 3$, it follows from \cite[Corollary 3.14]{VDKCentralExtensions} that $\Ker \pi=0$ unless $\ch k=3$ and $\tilde{\Phi}=A_2$. So this is the only case we have to consider. Since we are assuming $\rk\Phi\ge 2$, it follows that $u$ must be a bijection.
	
	In this last remaining case, in the notation of \cite{VDKCentralExtensions} we have $(\Ker \pi)_0=0$, which means that $\Ker\pi$ is generated by elements $Z_{\gamma}^*$, where $\gamma=\alpha+\beta$ is a degenerate sum for $\alpha$,$\beta\in\Phi=\tilde{\Phi}$. By relations from \cite[Theorem 3.5]{VDKCentralExtensions}, we have $Z_{\gamma}^*=\pm[X_{\alpha},X_{\beta}]$, where $\pi(X_{\alpha})\in\g_{\alpha}$ and $\pi(X_{\beta})\in\g_{\beta}$.
	
	Since $\alpha+\beta\notin\Phi$, applying R2 and Sublemma E we get
	$$
	\pi^\dag(Z_{\gamma}^*)\in [\Ker\tau+\g_{\alpha},\Ker\tau+\g_{\beta}]_{\dag}=0\tp
	$$

Therefore, we have $\Ker \pi\le\Ker\pi^\dag$.
	
\end{proof}

 \subsection*{''Diagonal'' automorphisms of classical Lie algebras}$\phan$

 \smallskip

 In the notation of the previous subsection consider the following functor from the category of rings to the category of groups.

 \begin{multline*}
 	\Diag(R)=\{(\sigma_{\alpha})_{\alpha\in\Phi}\colon \sigma_{\alpha}\in\Aut_R(\g_{\alpha}(R))\text{ and }\forall\alpha,\beta\in\Phi\; \forall x\in \g_{\alpha}(R)\; y\in\g_{\beta}(R)\\ (\alpha+\beta\in\Phi\;\wedge\; \alpha\notin\{-2\beta, -1/2\beta\})\Rightarrow \sigma_{\alpha+\beta}([x,y])=[\sigma_{\alpha}(x),\sigma_{\beta}(y)]\}\tc
 \end{multline*}

where $\Aut_R(\g_{\alpha}(R))$ is a group of automorphisms of $\g_{\alpha}(R)$ as an $R$-module.

It is easy to see that $\Diag(-)$ is an affine group scheme of finite type over $\Z$.

Let $G(\tilde{\Phi},-)$ be the adjoint Chevalley--Demazure group scheme of type $\tilde{\Phi}$; and let \mbox{$L(-)\le G(\tilde{\Phi},-)$} be the product of the subsystem subgroup for the subsystem $u^{-1}(0)$ with the torus. It is known that $L(-)$ is a smooth closed subgroup-subscheme of $G(\tilde{\Phi},-)$.

The adjoint action of $L(R)$ preserves the subspaces $\g_{\alpha}(R)$; hence it defines a group-scheme homomorphism $\map{\ad}{L(-)}{\Diag(-)}$.

{\lem\label{Clopen} Suppose that $\rk\Phi\ge 2$. Let $\mathcal{S}=\{2,3\}$ if $\tilde{\Phi}=G_2$; and $\mathcal{S}=\{2\}$ otherwise. Then the homomorphism $\map{\ad}{L_{\Z[\mathcal{S}^{-1}]}(-)}{\Diag_{\Z[\mathcal{S}^{-1}]}(-)}$ is simultaneously a closed embedding and an open embedding.}

\begin{proof}
	We will apply Lemma~\ref{Waterhouse}; let us verify that all its assumptions hold. Firstly, the scheme $L(-)$ is a semidirect product of Chevalley--Demazure scheme with a torus; hence, it is smooth; hence, it is flat. Secondly, for any $\Z[\mathcal{S}^{-1}]$-algebra $R$ the action of $G(\tilde{\Phi},R)$ on $\g(R)$ is faithful, and $\g(R)$ is generated as a Lie algebra by the root subspaces $\g_{\alpha}(R)$ for $\alpha\in\Phi$; hence, the map $\map{\ad_R}{L(R)}{\Diag(R)}$ is injective, which covers the assumption~(2).
	
	Now it remains to prove the assumption (1). Let $K$ be an algebraically closed field with $\ch K\notin\mathcal{S}$. It is easy to see that $\dim(L_K)=\rk\tilde{\Phi}+|u^{-1}(0)\cap \tilde{\Phi}|$. Let us now calculate $\dim_K\Lie(\Diag_K)$.
	
	Let us denote $\g=\g(K)$ and $\g_{\alpha}=\g_{\alpha}(K)$ for all $\alpha\in\Phi$. From the definition of the scheme $\Diag$ we obtain
	
	 \begin{multline*}
		\Lie(\Diag_K)=\{(\delta_{\alpha})_{\alpha\in\Phi}\colon \delta_{\alpha}\in\End_K(\g_{\alpha})\text{ and }\forall\alpha,\beta\in\Phi\; \forall x\in \g_{\alpha}\; y\in\g_{\beta}\\ (\alpha+\beta\in\Phi\;\wedge\; \alpha\notin\{-2\beta, -1/2\beta\})\Rightarrow \delta_{\alpha+\beta}([x,y])=[\delta_{\alpha}(x),y]+[x,\delta_{\beta}(y)]\}\tc
	\end{multline*}

	where $\End_K(\g_{\alpha})$ is the space of endomorphism of $\g_{\alpha}$ as a vector space.
	
	It follows from Lemma~\ref{LieCopresent} that any such collection $(\delta_{\alpha})_{\alpha\in\Phi}$ uniquely define a derivation $\delta\in\Der_K(\g)$ of the algebra $\g$.
	
	\medskip
	
	\underline{Sublemma} Let $\g_{\ad}$ be the tangent Lie algebra of the scheme $G(\tilde{\Phi},-)_K$. Then the action of $G(\tilde{\Phi},-)_K$ on $\g$ induce an isomorphism $\g_{\ad}\toiso\Der_K(\g)$.
	
	\medskip
	
	\underline{Proof of Sublemma.} The scheme $G(\tilde{\Phi},-)_K$ acts on $\g$ faithfully; hence the corresponding map $\g_{\ad}\to\Der_K(\g)$ is injective. Thus it remains to prove that $\dim_K\g_{\ad}=\dim_K\Der_K(\g)$.
	
	Since the simply connected and the adjoint Chevalley--Demazure schemes have the same dimension and are both smooth, it follows that $\dim_K\g_{\ad}=\dim_K \g$. Further, clearly, we have $\dim_K\Der_K(\g)=\dim_K \g-\dim_K Z(\g)+\dim_K H^1(\g,\g)$, where $\Z(\g)$ is the center of $\g$ and $H^1(\g,\g)$ is the space of outer derivations of $\g$. Thus it remains to prove that $\dim_K Z(\g)=\dim_K H^1(\g,\g)$.
	
	The algebra $\Z(\g)$ is the tangent algebra to the scheme-center of the simply connected Chevalley--Demazure scheme. Therefore, when this center is the scheme $\mu_n$ of $n$-th roots of unity, we have $\dim_K Z(\g)=1$ if $\ch K\mid n$ and $\dim_K Z(\g)=0$ otherwise. In the only remaining case, where $\Phi=D_{2l}$ and the scheme-center is $\mu_2\times\mu_2$ we have $\dim_K Z(\g)=2$ if $\ch K=2$ (which is not the case by assumption) and $\dim_K Z(\g)=0$ otherwise.
	
	The space $H^1(\g,\g)$ is calculated in \cite[Proposition 4.2]{IbarevFirstH} under the same assumptions on $\ch K$ that we have; and its dimension does indeed coincide with $\dim_K Z(\g)$.
	
	\underline{Sublemma is proved.}
	
	\medskip
	
	Sublemma together with what we show previously imply that $\Lie(\Diag_K)$ is isomorphic to the subalgebra of $\g_{\ad}$ that consists of element whose adjoint action preserves the subspaces~$\g_{\alpha}$. It is easy to see that this is exactly the subspace of elements whose components in root subalgebras with respect to grading by $\Phi$ vanish. Therefore, we have $\dim_K\Lie(\Diag_K)=\rk\tilde{\Phi}+|u^{-1}(0)\cap \tilde{\Phi}|$.

\end{proof}

\section{Over fields}
\label{FieldSec}

The following proposition combines the known results on abstract isomorphisms of isotropic simple groups over fields.

{\prop\label{BT} Let $G_1$ and $G_2$ be absolutely simple adjoint isotropic group schemes over fields $K_1$ and $K_2$. Let $\tilde{\Phi_1}$ and $\tilde{\Phi_2}$ be the corresponding absolute root systems. If the corresponding elementary groups are finite, then assume that they are simple {\rm (}i.e. they are not $E(A_1,\F_2)$, $E(A_1,\F_3)$, $E(B_2,\F_2)$, $E(G_2,\F_2)$ or $E(\phan^2 A_2,\F_2)${\rm )}. Then the only situations where isomorphism between $E_1(K_1)$ and $E_1(K_2)$ may not arise from a field isomorphism $\ph\colon K_1\toiso K_2$ and an $K_2$-group-scheme isomorphism $\Theta\colon \phan^\ph G_1\toiso G_2$ are the following:

$\bullet$ $\{\tilde{\Phi_1},\tilde{\Phi_2}\}=\{B_n,C_n\}$, $K_1\simeq K_2$, $\ch(K_1)=2$;

$\bullet$ $\tilde{\Phi_1}=\tilde{\Phi_2}=F_4$, $K_1\simeq K_2$, $\ch(K_1)=2$;

$\bullet$ $\tilde{\Phi_1}=\tilde{\Phi_2}=G_2$, $K_1\simeq K_2$, $\ch(K_1)=3$;

$\bullet$ $E(A_1,\F_4)\simeq E(A_1,\F_5)$;

$\bullet$ $E(A_1,\F_7)\simeq E(A_2,\F_2)$;

$\bullet$ $E(B_2,\F_3)\simeq E(\phan^2 A_3,\F_2)$.}

\begin{proof}
	Clearly, the fields $K_1$ and $K_2$ are either both finite, or both infinite. The case of the infinite fields is covered by \cite[Theorem 8.11]{BorelTitsHomo}. For the automorphisms of finite simple groups of Lie type see \cite[Theorem 3.2]{SteinbergAutFinite}, and for the isomorphisms between different groups see the list of finite simple groups.
	\end{proof}

\section{Corresponding maximal ideals and the proof of item (1) of Theorem~\ref{main}}
\label{MaximalSec}

In this section we are assuming the assumptions of Theorem~\ref{main}.

\smallskip

Notation $N^{(i)}_I=G_i(R_i,I)\cap E_i(R_i)$, $I\unlhd R_i$.

{\lem\label{MaxCorr} There is a one-to-one correspondence between the maximal ideals of $R_1$ and $R_2$, such that $\theta(N^{(1)}_I)=N^{(2)}_J$ whenever $I$ corresponds to $J$.}
\begin{proof}
	It follows from \cite[Theorem 1.1]{StavStep} that the subgroups $N^{(i)}_I$ for maximal ideals $I\unlhd R_i$ are precisely the maximal proper normal subgroups of $E_i(R_i)$. Clearly, isomorphism of groups maps maximal proper normal subgroups to maximal proper normal subgroups.
	\end{proof}

Now we prove item (1) of Theorem~\ref{main}.

\begin{proof}
Take maximal ideal $I\unlhd R_1$ and $J\unlhd R_2$ that correspond to each other as in Lemma~\ref{MaxCorr}. Then $ E_1(R_1/I)\simeq E_1(R_1)/N^{(1)}_I\simeq E_2(R_2)/N^{(2)}_J\simeq E_2(R_2/J)$. Then Proposition~\ref{BT} implies that $\tilde{\Phi_1}=\tilde{\Phi_2}$ unless we are at one of the exceptional cases. Suppose that $\tilde{\Phi_1}\ne\tilde{\Phi_2}$. The assumptions on isotropic rank, the assumption $2\in R_2^*$ and the assumption that $6\in R_2^*$ if $\tilde{\Phi_2}=G_2$ eliminate all the exceptions except the last one. Again since $2\in R_2^*$ we must have  $\tilde{\Phi_1}=A_3$, $\tilde{\Phi_2}=B_2$, $R_1/I\simeq \F_2$ and $R_2/J\simeq \F_3$. Finally, since we can take any pair of corresponding maximal ideals, we actually have $R_1/\M\simeq \F_2$ for all maximal ideals $\M\unlhd R_1$ and $R_2/\M\simeq \F_3$ for all maximal ideals $\M\unlhd R_2$.
\end{proof}

\section{Infinitesimal adjustment of unipotents}
\label{InfinitesimalSecI}

In this section the following {\bf umbrella assumptions} are in place:

$\bullet$ $\tilde{\Phi}$ is a reduced irreducible crystallographic root system; $\Phi$ is irreducible crystallographic root system of rank at least 2 (possibly of type $BC$);

$\bullet$ $\map{u}{\tilde{\Phi}\cup\{0\}}{\Phi\cup\{0\}}$ comes from a Tits index;

$\bullet$ $R$ and $S$ are commutative rings;

$\bullet$ $2\in S^*$;

$\bullet$ $I\unlhd S$ is an ideal such that $I^2=0$;

$\bullet$ $G$ is an absolutely simple adjoint group scheme over $R$ with absolute root system $\tilde{\Phi}$ that admits an isotropic pinning $(T^{(1)},\Phi)$ such that the corresponding map ${\tilde{\Phi}\sqcup\{0\}}\to{\Phi\sqcup\{0\}}$ is $u$, and this pinning satisfies the assumption (a)--(d) in the second bullet in the statement of Theorem~\ref{main};

$\bullet$ $E(R)$ --- the elementary subgroup of $G(R)$;

$\bullet$ $G_{\alpha}^{(1)}=G_{\alpha}^{(1)}(R)$ for $\alpha\in\Phi$ is the corresponding root subgroup,  $h_{\alpha}^{(1)}$ are the elements from Lemma~\ref{haExist}; note that since we are assuming the square formula, we have $h_{\alpha}^{(1)}\in E(R)$;

$\bullet$ $G(\tilde{\Phi},-)$ is an adjoint Chevalley--Demazure group scheme of type $\tilde{\Phi}$.

$\bullet$ $(T^{(2)},\Phi)$ is an isotropic pre-pinning on $G(\tilde{\Phi},-)$ that is contained in the pinning and such that the corresponding map ${\tilde{\Phi}\sqcup\{0\}}\to{\Phi\sqcup\{0\}}$ is $u$. Note that by Lemma~\ref{WeylForSplit} this isotropic pre-pinning is actually an isotropic pinning.

$\bullet$ $G_{\alpha}^{(2)}$ for $\alpha\in\Phi$ is the corresponding root subgroup-subscheme (meaning that we have $G_{\alpha}^{(2)}=\prod_{u(\beta)\in\{\alpha,2\alpha\}} U_{\beta}$), and $h_{\alpha}^{(2)}$ are the elements from Lemma~\ref{haExist};

$\bullet$ $L^{(2)}(-)\le G(\tilde{\Phi},-)$ is the product of the subsystem subgroup for the subsystem $u^{-1}(0)$ with the torus.

$\bullet$  $C^{\mathrm{us}}=\bigcap_{2\alpha\notin\Phi} C_{L^{(2)}}(G_{\alpha}^{(2)})$ is the scheme centralizer of non-ultrashort root subgroups (so it is only nontrivial if $\Phi$ is of type $BC$).

$\bullet$ $\map{\rho_I}{G(\tilde{\Phi},S)}{G(\tilde{\Phi},S/I)}$ is the reduction homomorphism;

$\bullet$ $\map{\theta}{E(R)}{G(\tilde{\Phi},S)}$ is a homomorphism;

$\bullet$ $\rho_I(\theta(h^{(1)}_{\alpha}))=\rho_I(h^{(2)}_\alpha)$ for all $\alpha\in\Phi$;

$\bullet$  $\rho_I(\theta(G^{(1)}_{\alpha}))\sub G^{(2)}_{\alpha}(S/I)$ for all $\alpha\in\Phi$; and $\rho_I(\theta(G^{(1)}_{\alpha}))$ generate $G^{(2)}_{\alpha}(S/I)$ (or $G^{(2)}_{\alpha}(S/I)/G^{(2)}_{2\alpha}(S/I)$ if $\alpha$ is ultrashort) as an $S/I$-module.

\medskip

Since $I^2=0$, it follows that $G(\tilde{\Phi},S,I)$ is an abelian group with the structure of an $S$-module. We will constantly use this fact without mentioning. Any element $\eps\in G(\tilde{\Phi},S,I)$ can be uniquely decomposed into a product of elements from $L^{(2)}(S)\cap G(\tilde{\Phi},S,I)$ and from $G^{(2)}_{\beta}(S)\cap G(\tilde{\Phi},S,I)$ for all $\beta\in \Phi\sm 2\Phi$. Since $G(\tilde{\Phi},S,I)$ is abelian, the factors of this decomposition are well defined independently on the order of roots. We will call these factors the components of $\eps$.

{\lem\label{GiantLemma} Under the umbrella assumptions of this Section, assume additionally that $\theta(h^{(1)}_{\alpha})=h^{(2)}_\alpha$ for all non-ultrashort $\alpha\in\Phi$. Then
	\begin{enumerate}
		\item if $\Phi$ is not of type $C$ or $BC$ (including $C_2=B_2$), then  we have $\theta(G^{(1)}_{\alpha})\sub G^{(2)}_{\alpha}(S)$ for all $\alpha\in\Phi$;
		
		\item if $\Phi$ is of type $C$ or $BC$, then there exist an element
		$$
		g\in \prod_{\alpha \text{ is long}} \left(G^{(2)}_{\alpha}(S)\cap G(\tilde{\Phi},S,I)\right)
		$$
		such that for the homomorphism $\theta'(x)=\theta(x)^g$ we have $\theta'(G^{(1)}_{\alpha})\sub G^{(2)}_{\alpha}(S)$ for all $\alpha\in\Phi$;
	\end{enumerate}
	 }
\begin{proof}
	Clearly $\map{\rho_I}{G^{(2)}_{\alpha}(S)}{G^{(2)}_{\alpha}(S/I)}$ is surjective. Thus for any $x\in G_{\alpha}^{(1)}$ we have $\theta(x)=y_x\eps_x$, where $y_x\in G^{(2)}_{\alpha}(S)$ and $\eps_x\in G(\tilde{\Phi},S,I)$. Further we can decompose $\eps_x$ as $\eps_x=\prod_{\beta\in\Phi\sqcup\{0\}\sm 2\Phi} \eps_{x,\beta}$, where $\eps_{x,0}\in L^{(2)}(S)\cap G(\tilde{\Phi},S,I)$ and $\eps_{x,\beta}\in G^{(2)}_{\beta}(S)\cap G(\tilde{\Phi},S,I)$ for all $\beta\in \Phi\sm 2\Phi$.
	
	\medskip
	
	\underline{Sublemma A.} Let $x\in G_{\alpha}^{(1)}$ and $\beta\in \Phi\sm 2\Phi$. Suppose that for some non-ultrashort root $\gamma\in\Phi$ we have $\<\alpha,\gamma\>$ is even and $\<\beta,\gamma\>$ is odd. Then for non-ultrashort $\beta$ we have $\eps_{x,\beta}=e$; and if for ultrashort $\beta$ we have $\eps_{x,\beta}\in G^{(2)}_{2\beta}(S)$.
	
	\medskip
	
	\underline{Proof of Sublemma A.} Since $\<\alpha,\gamma\>$ is even we have $x^{h^{(1)}_{\gamma}}=x$, hence $\theta(x)^{h^{(2)}_{\gamma}}=\theta(x)$. On the other hand, we have $\theta(x)^{h^{(2)}_{\gamma}}=y_x\prod_{\beta'\in\Phi\sqcup\{0\}\sm 2\Phi} (\eps_{x,\beta'})^{h^{(2)}_{\gamma}}$. Thus we have $(\eps_{x,\beta'})^{h^{(2)}_{\gamma}}=\eps_{x,\beta'}$ for every $\beta'\in\Phi\sqcup\{0\}\sm 2\Phi$. In particular, we have $\eps_{x,\beta}=(\eps_{x,\beta})^{h^{(2)}_{\gamma}}=\eps_{x,\beta}\cdot (-1)$. If $\beta$ is not ultrashort, then we have $\eps_{x,\beta}\cdot 2=e$, and since $2\in S^*$, we have  $\eps_{x,\beta}=e$. If $\beta$ is ultrashort, then the similar argument shows that the image of $\eps_{x,\beta}$ in $G^{(2)}_{\beta}(S)/G^{(2)}_{2\beta}(S)$ is trivial, i.e. we have $\eps_{x,\beta}\in G^{(2)}_{2\beta}(S)$. \underline{Sublemma A is proved.}
	
	\medskip
	
	\underline{Sublemma B.} Let $x\in G_{\alpha}^{(1)}$, where $\alpha$ is not ultrashort, and $\beta\in \Phi\sm \R\alpha$ is also not ultrashort. Suppose that for some non-ultrashort root $\gamma\in\Phi$ we have $\<\alpha,\gamma\>$ is odd and $\<\beta,\gamma\>$ is even. If $\beta\in 2\Phi$, then assume that $\eps_{x,\beta/2}\in G_{\beta}^{(2)}(S)$ and denote $\eps_{x,\beta}=\eps_{x,\beta/2}$. Assume additionally that $\eps_{x,\beta'}=e$ for all $\beta'\in\Phi$ such that either $\beta\in \beta'+\N_+\alpha$ or $\beta/2\in \beta'+\N_+\alpha$. Then we have $\eps_{x,\beta}=e$.
	
	\medskip
	
	\underline{Proof of Sublemma B.} Since $\<\alpha,\gamma\>$ is odd and $\alpha$ is not ultrashort we have $x^{h^{(1)}_{\gamma}}=x^{-1}$, hence
	$$
	\theta(x)^{h^{(2)}_{\gamma}}=\theta(x)^{-1}=\left(\prod_{\beta'\in\Phi\sqcup\{0\}\sm 2\Phi} (\eps_{x,\beta'})^{-1}\right)y_x^{-1}=y_x^{-1}\!\!\prod_{\beta'\in\Phi\sqcup\{0\}\sm 2\Phi} (\eps_{x,\beta'})^{-1}\cdot\!\!\prod_{\beta'\in\Phi\sqcup\{0\}\sm 2\Phi} [\eps_{x,\beta'},y_x]\tp
	$$
	 On the other hand, we have $\theta(x)^{h^{(2)}_{\gamma}}=y_x^{-1}\prod_{\beta'\in\Phi\sqcup\{0\}\sm 2\Phi} (\eps_{x,\beta'})^{h^{(2)}_{\gamma}}$. Note that we have $[\eps_{x,\beta'},y_x]\in\prod_{\beta''\in \beta'+\N_+\alpha} G^{(2)}_{\beta''}(S)$. Therefore, we obtain $(\eps_{x,\beta})^{h^{(2)}_{\gamma}}=(\eps_{x,\beta})^{-1}$. Since $\<\beta,\gamma\>$ is even, we also have $(\eps_{x,\beta})^{h^{(2)}_{\gamma}}=\eps_{x,\beta}$. Thus $\eps_{x,\beta}=(\eps_{x,\beta})^{-1}$, which since $2\in S^*$ implies $\eps_{x,\beta}=e$. \underline{Sublemma B is proved.}
	
	\medskip
	
	\underline{Local Definition.} We say that $\alpha$,$\beta\in\Phi$ are {\it parity-equivalent} if for any non-ultrashort $\gamma\in\Phi$ we have $\<\alpha,\gamma\>\equiv\<\beta,\gamma\>\mod 2$.
	
	\medskip
	
	\underline{Sublemma C.} Any two parity-equivalent roots are either collinear or orthogonal.
	
	\medskip
	
	\underline{Proof of Sublemma C} is by direct case-by-case computation.
	
	\medskip
	
	\underline{Sublemma D.} Let $x\in G_{\alpha}^{(1)}$ and $x'\in G_{\beta}^{(1)}$, where $\alpha+\beta\notin \Phi\cup\{0\}$. Then we have $[y_x^{-1},\eps_{x'}^{-1}]=[y_{x'}^{-1},\eps_x^{-1}]$.
	
	\medskip
	
	\underline{Proof of Sublemma D.}
	\begin{multline*}
	e=\theta(e)=\theta([x^{-1},(x')^{-1}])=[\eps_x^{-1}y_x^{-1},\eps_{x'}^{-1}y_{x'}^{-1}]=[\eps_x^{-1}y_x^{-1},\eps_{x'}^{-1}]\cdot[\eps_x^{-1}y_x^{-1},y_{x'}^{-1}]^{\eps_x'}=\\=[\eps_x^{-1}y_x^{-1},\eps_{x'}^{-1}]\cdot[\eps_x^{-1}y_x^{-1},y_{x'}^{-1}]=[y_x^{-1},\eps_{x'}^{-1}]^{\eps_x}\cdot[\eps_x^{-1},\eps_{x'}^{-1}]\cdot [y_x^{-1},y_{x'}^{-1}]^{\eps_x}\cdot[\eps_x^{-1},y_{x'}^{-1}]=\\=[y_x^{-1},\eps_{x'}^{-1}]\cdot [\eps_x^{-1},y_{x'}^{-1}]\tp
\end{multline*}
\underline{Sublemma D is proved.}

	\medskip

	\underline{Notation}: for every non-ultrashort root $\alpha\in\Phi$ we fix elements $a_{\alpha}$,$c_{\alpha}\in G_{\alpha}^{(1)}$ and $b_{\alpha}\in G_{-\alpha}^{(1)}$ such that $w_{\alpha}^{(1)}=a_{\alpha}b_{\alpha}c_{\alpha}$ is a Weyl element. Then we set $\tilde{w}_{\alpha}=y_{a_{\alpha}}y_{b_{\alpha}}y_{c_{\alpha}}$.

	\medskip

	\underline{Sublemma E.} For any $\alpha$,$\beta\in\Phi$, where $\alpha$ is non-ultrashort we have
	$$
	\leftact{\tilde{w}_{\alpha}}{(G_{\beta}^{(2)}(S))}\sub G_{s_{\alpha}(\beta)}^{(2)}(S)G(\tilde{\Phi},S,I)
	$$
	and
	$$
	(G_{\beta}^{(2)}(S))^{\tilde{w}_{\alpha}}\sub G_{s_{\alpha}(\beta)}^{(2)}(S)G(\tilde{\Phi},S,I)\tp
	$$

\medskip

	\underline{Proof of Sublemma E.}
	
	We write the proof for the first statement. The second statement is similar.
	
	We must prove that $\rho_I(\leftact{\tilde{w}_{\alpha}}{(G_{\beta}^{(2)}(S))})\sub G_{s_{\alpha}(\beta)}^{(2)}(S/I)$. We can do it by showing that elements from $\rho_I(\leftact{\tilde{w}_{\alpha}}{(G_{\beta}^{(2)}(S))})$ satisfy the formulas form Lemma~\ref{EgorFormulas} taking $X_{\gamma}=\rho_I(\theta(G_{\gamma}^{(1)}))$ for any $\gamma\in\Phi$.
	
	Let us check the first formula suppose that $\gamma\in\Phi$ is such that $s_\alpha(\beta)+\gamma\notin\Phi\cup\{0\}$, let $g\in G_{\beta}^{(2)}(S)$, and $x_{\gamma}\in G_{\gamma}^{(1)}$. Then we have
	$$
	[\rho_I(\leftact{\tilde{w}_{\alpha}}{g}),\rho_I(\theta(x_\gamma))]=\leftact{\rho_I(\tilde{w}_{\alpha})}{[\rho_I(g),\rho_I(\theta(x_{\gamma}^{w_{\alpha}^{(1)}}))]}\le \leftact{\rho_I(\tilde{w}_{\alpha})}{[G_{\beta}^{(2)}(S/I),G_{s_\alpha(\gamma)}^{(2)}(S/I)]}=e\tp
	$$
	
	The other two formulas can be verified similarly.

	\underline{Sublemma E is proved.}
	
	\medskip
	\underline{Sublemma F.} Let $\alpha$ be a long root, and $x\in G_{\alpha}^{(1)}$. Then we have $\theta(x)=y_x\eps_{x,\alpha}\eps_{x,0}$ (here for $\alpha\in 2\Phi$ we denoted $\eps_{x,\alpha}=\eps_{x,\alpha/2}$ and we state that $\eps_{x,\alpha}\in G_{\alpha}^{(2)}(S)$). Moreover, if $\Phi$ is not of type $C$ or $BC$, then we have $\eps_{x,0}=e$; and if $\Phi$ is of type $C$ or $BC$, then $\eps_{x,0}$ commutes with $G_{\beta}^{(2)}(S)$ for all non-ultrashort roots $\beta\in\Phi$ such that $\beta\perp\alpha$.
	
	\medskip
	
	\underline{Proof of Sublemma F.} It is easy to check case by case that if $\alpha$ is a long root, then for any $\beta\in\Phi$ that is not parity-equivalent to $\alpha$ there exist a non-ultrashort root $\gamma\in\Phi$ such that $\<\alpha,\gamma\>$ is even and $\<\beta,\gamma\>$ is odd. Therefore, Sublemma A implies that
	$$
	\theta(x)=y_x\eps_{x,0}\prod_{\beta\sim\alpha}\eps_{x,\beta}\tc
	$$
	where $\sim$ stands for parity-equivalence, and $\eps_{x,\beta}\in G_{\beta}^{(2)}(S)\cap G(\tilde{\Phi},S,I)$ (here for $\beta\in 2\Phi$ we denoted $\eps_{x,\beta}=\eps_{x,\beta/2}$ and Sublemma A implies that $\eps_{x,\beta}\in G_{\beta}^{(2)}(S)$).
	
	Note that by Sublemma C all the roots in the parity-equivalence class of $\alpha$ are either $\pm\alpha$ or orthogonal to $\alpha$.
	
	The next step is to show that $\eps_{x,-\alpha}=e$. Choose $\beta\in\Phi$ as follows: if $\alpha$ lies in a subsystem of type $A_2$, then let $\beta$ be a long root with $\angle(\alpha,\beta)=\pi/3$; otherwise let $\beta$ be a short root with $\angle(\alpha,\beta)=\pi/4$. By Sublemma D for any $x'\in G_{\beta}^{(1)}$ we have $[y_x^{-1},\eps_{x'}^{-1}]=[y_{x'}^{-1},\eps_x^{-1}]$. Clearly, the component in $G_{\beta-\alpha}^{(2)}(S)$ at the left hand side vanishes; and at the right hand side such component comes from $[y_{x'}^{-1},\eps_{x,-\alpha}^{-1}]$. By umbrella assumptions $\rho_I(\theta(G^{(1)}_{\beta}))$ generate $G^{(2)}_{\beta}(S/I)$ as an $S/I$-module; hence, all possible $y_{x'}$ generate $G^{(2)}_{\beta}(S)$ as an $S$-module; hence $[\eps_{x,-\alpha},G^{(2)}_{\beta}(S)]$ is either trivial or lies in $G^{(2)}_{2\beta-\alpha}(S)$ if $\beta$ is short. The commutator (taken modulo $G^{(2)}_{2\beta-\alpha}(S)$ if $\beta$ is short) is a non-degenerate pairing of $G^{(2)}_{-\alpha}(S)$ and $G^{(2)}_{\beta}(S)$ (if the root $\beta$ is long, then this is because of the existence of Weyl elements, otherwise it is by \cite[Lemma~3]{VoronDiophantine}), thus we have $\eps_{x,-\alpha}=e$.
	
	The next step is to prove what is claimed about $\eps_{x,0}$. First let us consider the case, where $\Phi$ is not of type $C$ or $BC$. In this case there is a non-ultrashort root $\gamma\in\Phi$ such that $\<\alpha,\gamma\>$ is odd. Therefore, we have
	$$
	\theta(x)^{-1}=\theta(x^{-1})=\theta(x^{h_{\gamma}^{(1)}})=\theta(x)^{h_{\gamma}^{(2)}}\tp
	$$
	Now using what we prove before, we can compare the components in $L^{(2)}(S)$ on the left hand side and on the right hand side of this equality and conclude that $\eps_{x,0}=\eps_{x,0}^{-1}$; hence, since $2\in S^*$ we have $\eps_{x,0}=e$.
	
	Now suppose that $\Phi$ is of type $C$ or $BC$. We must prove that $\eps_{x,0}$ commutes with $G_{\beta}^{(2)}(S)$ for all non-ultrashort roots $\beta\in\Phi$ such that $\beta\perp\alpha$. Since all possible $y_{x'}$, where $x'\in G^{(1)}_{\beta}$, generate $G^{(2)}_{\beta}(S)$ as an $S$-module, it is enough to show that $\eps_{x,0}$ commutes with all such~$y_{x'}$.
	
	First consider the case, where $\beta$ is short. By Sublemma D we have $[y_x^{-1},\eps_{x'}^{-1}]=[y_{x'}^{-1},\eps_x^{-1}]$. It follows from what we proved before that the commutator $[y_{x'}^{-1},\eps_{x,0}^{-1}]$ is equal to the component in $G_{\beta}^{(2)}(S)$ at the right hand side of this equality; and at the left hand side such component is trivial, since $\beta-\alpha\notin\Phi\cup\{0\}$. Therefore, $y_{x'}^{-1}$ commutes with $\eps_{x,0}^{-1}$; hence $y_{x'}$ commutes with $\eps_{x,0}$.
	
	Now suppose that $\beta$ is a long root. By Sublemma D we have $[y_x^{-1},\eps_{x'}^{-1}]=[y_{x'}^{-1},\eps_x^{-1}]$. Let us compare the components in $G_{\beta}^{(2)}(S)$ at the right and left hand sides of this equality. At the left hand side such component is trivial, since $\beta-\alpha\notin\Phi\cup\{0\}$. At the right hand side such component is equal to $u(x')v(x')$, where $u(x')=[y_{x'}^{-1},\eps_{x,0}^{-1}]$ and $v(x')$ is the component in $G_{\beta}^{(2)}(S)$ of the element $[y_{x'}^{-1},\eps_{x,-\beta}^{-1}]$. Therefore, we have $u(x')v(x')=e$. Similarly, we have $u((x')^{-1})v((x')^{-1})=e$. Without loss of generality, we may assume that $y_{(x')^{-1}}=y_{x'}^{-1}$. In this case, it is easy to see that $u((x')^{-1})=u(x')^{-1}$ and $v((x')^{-1})=v(x')$; hence, we have $u(x')^{-1}v(x')=e$. From this equality and the equality $u(x')v(x')=e$, we get $u(x')^2=e$; hence, since $2\in S^*$ we have $u(x')=e$. Therefore, $y_{x'}^{-1}$ commutes with $\eps_{x,0}^{-1}$; hence $y_{x'}$ commutes with $\eps_{x,0}$.
	
	The final step is to show that $\eps_{x,\beta}=e$, where $\beta\in\Phi$ is such that $\beta\sim\alpha$ and $\beta\perp\alpha$. By Sublemma D we have $[y_x^{-1},\eps_{b_{\beta}}^{-1}]=[y_{b_{\beta}}^{-1},\eps_x^{-1}]$. Let us compare the components in $G_{-\beta}^{(2)}(S)$ in this equality. At the left hand side such component is trivial, since $\beta-\alpha\notin\Phi\cup\{0\}$. At the right hand side contribution to such component may come only from $[y_{b_{\beta}}^{-1},\eps_{x,0}^{-1}]$, which is trivial as we prove previously, and from $[y_{b_{\beta}}^{-1},\eps_{x,\beta}^{-1}]$; hence, we conclude that $[y_{b_{\beta}}^{-1},\eps_{x,\beta}^{-1}]$ has trivial component in $G_{-\beta}^{(2)}(S)$. Now let us prove that $(\eps_{x,\beta}^{-1})^{\tilde{w}_{\beta}}\in G_{-\beta}^{(2)}(S)$. Indeed, since $\beta$ is neither ultrashort nor a short root when $\Phi$ is of type $C$ or $BC$, by Lemma~\ref{EgorFormulas} it is enough to show that $(\eps_{x,\beta}^{-1})^{\tilde{w}_{\beta}}$ commutes with $G_{\gamma}^{(2)}(S)$ for all $\gamma$ such that $\gamma-\beta\notin\Phi\cup\{0\}$. It is the same that to show that $\eps_{x,\beta}^{-1}$ commutes with $\leftact{\tilde{w}_{\beta}}{(G_{\gamma}^{(2)}(S))}$. By Sublemma E we have $\leftact{\tilde{w}_{\beta}}{(G_{\gamma}^{(2)}(S))}\sub G_{s_{\beta}(\gamma)}^{(2)}(S)G(\tilde{\Phi},S,I)$. Since $\gamma-\beta\notin\Phi\cup\{0\}$, we have $s_{\beta}(\gamma)+\beta\notin\Phi\cup\{0\}$; hence $\eps_{x,\beta}^{-1}$ commutes with $G_{s_{\beta}(\gamma)}^{(2)}(S)$; and since $G(\tilde{\Phi},S,I)$ is abelian, $\eps_{x,\beta}^{-1}$ also commutes with $G(\tilde{\Phi},S,I)$. Therefore, we conclude that $(\eps_{x,\beta}^{-1})^{\tilde{w}_{\beta}}\in G_{-\beta}^{(2)}(S)$. Now it is easy to see that $(\eps_{x,\beta}^{-1})^{\tilde{w}_{\beta}}$ must be equal to the component of $[y_{b_{\beta}}^{-1},\eps_{x,\beta}^{-1}]$ in $G_{-\beta}^{(2)}(S)$, which we prove to be trivial. Therefore, we have $(\eps_{x,\beta}^{-1})^{\tilde{w}_{\beta}}=e$; hence, $\eps_{x,\beta}=e$.

	 \underline{Sublemma F is proved.}
	
	 \medskip
	
	 \underline{Sublemma G.} Suppose that $\Phi$ is of type $C$ or $BC$ and $\alpha$,$\beta\in\Phi$ be such that $\alpha$ is long $\beta$ is short and $\angle(\alpha,\beta)=\pi/4$. Then for any $x\in G^{(1)}_{\alpha-\beta}$ we have $\eps_{x,-\alpha}=\eps_{x,\beta-\alpha}=\eps_{x,0}=e$ (here if $\Phi$ is of type $BC$, then we denote $\eps_{x,-\alpha}=\eps_{x,-\alpha/2}$, note that by Sublemma A we have $\eps_{x,-\alpha}\in G_{-\alpha}^{(2)}(S)$).
	
	 \medskip
	
	 \underline{Proof of Sublemma G.}
	
	 First $\eps_{x,-\alpha}=e$ by Sublemma B. Now by Sublemma D for every $x'\in G^{(1)}_{\alpha-2\beta}$ we have $[y_x^{-1},\eps_{x'}^{-1}]=[y_{x'}^{-1},\eps_x^{-1}]$. Comparing components in $G^{(2)}_{-\beta}(S)$, while keeping in mind that by Sublemma F we have $\eps_{x'}=\eps_{x',\alpha-2\beta}\eps_{x',0}$, we get $[y_{x'}^{-1},\eps_{x,\beta-\alpha}^{-1}]\in G^{(2)}_{-\alpha}(S)$. By \cite[Lemma~3]{VoronDiophantine} the commutator taken modulo $G^{(2)}_{-\alpha}(S)$ is a non-degenerate paring between $G^{(2)}_{\beta-\alpha}(S)$ and $G^{(2)}_{\alpha-2\beta}(S)$; and since all possible $y_{x'}$ generate $G^{(2)}_{\alpha-2\beta}(S)$ as an $S$-module, we have $\eps_{x,\beta-\alpha}=e$.
	
	 Further we have
	 $$
	 \theta(x)^{-1}=\theta(x^{-1})=\theta(x^{h_{\alpha}^{(1)}})=\theta(x)^{h_{\alpha}^{(2)}}\tc
	 $$
	 so using what we prove before, we can compare the components in $L^{(2)}(S)$ on the left hand side and on the right hand side of this equality and conclude that $\eps_{x,0}=\eps_{x,0}^{-1}$; hence, since $2\in S^*$ we have $\eps_{x,0}=e$.
	
	 \underline{Sublemma G is proved.}
	
	\medskip
	
		\underline{Sublemma H.} Suppose that $\Phi$ is of type $C$ or $BC$, and let $\alpha\in\Phi$ be a long root. Then there exists an element $g_{-\alpha}\in G_{-\alpha}^{(2)}(S)\cap G(\tilde{\Phi},S,I)$ such that after we replace the homomorphism $\theta(-)$ with the homomorphism $\theta(-)^{g_{-\alpha}}$ we will have $\eps_{x,0}\in C^{\mathrm{us}}$ for every $x\in G_{\alpha}^{(1)}$ (if $\Phi$ is of type $C$ this means that $\eps_{x,0}=e$).
		
	\medskip
	
	\underline{Proof of Sublemma H.}
	
	Fix a short root $\beta\in\Phi$ such that $\angle(\alpha,\beta)=\pi/4$. Using the assumption (d) in the second bullet in the statement of Theorem~\ref{main}, we may assume without loss of generality that, if $\map{u}{\tilde{\Phi}\cup\{0\}}{\Phi\cup\{0\}}$ is not a bijection, then $b_{-\beta}$ commutes with $b_{\beta-\alpha}$.
	
	Note that $\leftact{\tilde{w}_{-\beta}}{(G^{(2)}_{2\beta-\alpha}(S)\cap G(\tilde{\Phi},S,I))}\sub G^{(2)}_{-\alpha}(S)\cap G(\tilde{\Phi},S,I)$. Indeed, by Lemma~\ref{EgorFormulas} it it enough to show that $\leftact{\tilde{w}_{-\beta}}{(G^{(2)}_{2\beta-\alpha}(S)\cap G(\tilde{\Phi},S,I))}$
	commutes with $G_{\gamma}^{(2)}(S)$ for all $\gamma$ such that  $\gamma-\alpha\notin\Phi\cup\{0\}$. It is the same that to show that $G^{(2)}_{2\beta-\alpha}(S)\cap G(\tilde{\Phi},S,I)$ commutes with $(G_{\gamma}^{(2)}(S))^{\tilde{w}_{-\beta}}$. By Sublemma E we have $(G_{\gamma}^{(2)}(S))^{\tilde{w}_{-\beta}}\sub G_{s_{\beta}(\gamma)}^{(2)}(S)G(\tilde{\Phi},S,I)$. Since $\gamma-\alpha\notin\Phi\cup\{0\}$, we have $s_{\beta}(\gamma)+2\beta-\alpha\notin\Phi\cup\{0\}$; hence $G^{(2)}_{2\beta-\alpha}(S)\cap G(\tilde{\Phi},S,I)$ commutes with $G_{s_{\beta}(\gamma)}^{(2)}(S)$; and since $G(\tilde{\Phi},S,I)$ is abelian it also commutes with $G(\tilde{\Phi},S,I)$.
	
	Hence, we have $G^{(2)}_{2\beta-\alpha}(S)\cap G(\tilde{\Phi},S,I)\sub (G^{(2)}_{-\alpha}(S)\cap G(\tilde{\Phi},S,I))^{\tilde{w}_{-\beta}}$.
	
	 Now for any $g_{-\alpha}\in G_{-\alpha}^{(2)}(S)\cap G(\tilde{\Phi},S,I)$ such that $g_{-\alpha}^{\tilde{w}_{-\beta}}\in G^{(2)}_{2\beta-\alpha}(S)\cap G(\tilde{\Phi},S,I)$ we have
	$$
	\theta(b_{-\beta})^{g_{-\alpha}}=y_{b_{-\beta}}^{g_{-\alpha}}\eps_{b_{-\beta}}=y_{b_{\beta}}\eps'_{b_{-\beta}}\tc
	$$
	where $\eps'_{b_{-\beta}}\in G(\tilde{\Phi},S,I)$ has component in $G^{(2)}_{2\beta-\alpha}(S)$ equal to $(g_{-\alpha}^{-1})^{\tilde{w}_{-\beta}}\eps_{b_{-\beta},2\beta-\alpha}$. Since $G^{(2)}_{2\beta-\alpha}(S)\cap G(\tilde{\Phi},S,I)\sub (G^{(2)}_{-\alpha}(S)\cap G(\tilde{\Phi},S,I))^{\tilde{w}_{-\beta}}$, it follows that $g_{-\alpha}$ can be chosen so that the last expression is trivial. That will be our $g_{-\alpha}$.
	
	From this moment suppose that we replaced the homomorphism $\theta(-)$ with the homomorphism $\theta(-)^{g_{-\alpha}}$; so will have $\eps_{b_{-\beta},2\beta-\alpha}=e$. Note that $\theta(h^{(1)}_{\gamma})$ has not changed for all $\gamma\in \Phi$; thus results of all the previous sublemmas still hold.
	
	\smallskip
	
	\underline{Subsublemma} Let $x\in G^{(1)}_{\alpha-\beta}$. Suppose that either $x$ commutes with $b_{-\beta}$, or $\eps_{b_{-\beta},\beta-\alpha}=e$. Then we have $\eps_{x,-\beta}=\eps_{x,\alpha-2\beta}=e$.
	
	\smallskip
	
	\underline{Proof of Subsublemma.}
	
	Since $[x^{-1},b_{-\beta}^{-1}]\in G^{(1)}_{\alpha}$ by Sublemma F we have $\theta([x^{-1},b_{-\beta}^{-1}])\in L^{(2)}(S)G^{(2)}_{\alpha}(S)$. On the other hand we have
	\begin{multline*}
		\theta([x^{-1},b_{-\beta}^{-1}])=[\theta(x)^{-1},\theta(b_{-\beta})^{-1}]=[\eps_x^{-1}y_x^{-1},\eps_{b_{-\beta}}^{-1}y_{b_{-\beta}}^{-1}]=[\eps_x^{-1}y_x^{-1},\eps_{b_{-\beta}}^{-1}]\cdot [\eps_x^{-1}y_x^{-1},y_{b_{-\beta}}^{-1}]^{\eps_{b_{-\beta}}}=\\=[y_x^{-1},\eps_{b_{-\beta}}^{-1}]^{\eps_x}\cdot [\eps_x^{-1},\eps_{b_{-\beta}}^{-1}] \cdot[y_x^{-1},y_{b_{-\beta}}^{-1}]^{\eps_x \eps_{b_{-\beta}}}\cdot [\eps_x^{-1},y_{b_{-\beta}}^{-1}]^{\eps_{b_{-\beta}}}=[y_x^{-1},\eps_{b_{-\beta}}^{-1}]\cdot [y_x^{-1},y_{b_{-\beta}}^{-1}]^{\eps_x \eps_{b_{-\beta}}}\cdot\\ \cdot [\eps_x^{-1},y_{b_{-\beta}}^{-1}]= [y_x^{-1},\eps_{b_{-\beta}}^{-1}]\cdot [\eps_x^{-1},y_{b_{-\beta}}^{-1}]\cdot [y_x^{-1},y_{b_{-\beta}}^{-1}]^{\eps_x \eps_{b_{-\beta}}[\eps_x^{-1},y_{b_{-\beta}}^{-1}]}=[y_x^{-1},\eps_{b_{-\beta}}^{-1}]\cdot [\eps_x^{-1},y_{b_{-\beta}}^{-1}]\cdot\\ \cdot [\eps_x^{-1}\eps_{b_{-\beta}}^{-1}[y_{b_{-\beta}}^{-1},\eps_x^{-1}],[y_x^{-1},y_{b_{-\beta}}^{-1}]]\cdot [y_x^{-1},y_{b_{-\beta}}^{-1}]\tp
	\end{multline*}
	Thus we obtain that $[y_x^{-1},\eps_{b_{-\beta}}^{-1}]\cdot [\eps_x^{-1},y_{b_{-\beta}}^{-1}]\cdot[\eps_x^{-1}\eps_{b_{-\beta}}^{-1}[y_{b_{-\beta}}^{-1},\eps_x^{-1}],[y_x^{-1},y_{b_{-\beta}}^{-1}]]\in L^{(2)}(S)G^{(2)}_{\alpha}(S)$. In particular its component in $G^{(2)}_{\beta}(S)$ is trivial. Let us show that the first and the third factors have trivial components in $G^{(2)}_{\beta}(S)$ so that we could conclude the same for the second factor. The first factor has trivial components in $G^{(2)}_{\beta}(S)$, because $\eps_{b_{-\beta},2\beta-\alpha}=e$. For the third factor we use our assumption that either $x$ commutes with $b_{-\beta}$, or $\eps_{b_{-\beta},\beta-\alpha}=e$. In case where $x$ commutes with $b_{-\beta}$ we have
	$$
	[y_x^{-1},y_{b_{-\beta}}^{-1}]\in [\theta(x)^{-1},\theta(b_{-\beta})^{-1}] G(\tilde{\Phi},S,I)=G(\tilde{\Phi},S,I)\tp
	$$
	Hence, since $G(\tilde{\Phi},S,I)$ is abelian the whole third factor is trivial. Now suppose that $\eps_{b_{-\beta},\beta-\alpha}=e$. The component of the third factor in $G^{(2)}_{\beta}(S)$ depends on the component of $\eps_x^{-1}\eps_{b_{-\beta}}^{-1}[y_{b_{-\beta}}^{-1},\eps_x^{-1}]$ in $G^{(2)}_{\beta-\alpha}(S)$. By assumption we have $\eps_{b_{-\beta},\beta-\alpha}=e$; by Sublemma G we have $\eps_{x,\beta-\alpha}=e$; and $[y_{b_{-\beta}}^{-1},\eps_x^{-1}]$ has trivial component in $G^{(2)}_{\beta-\alpha}(S)$ because by Sublemma G we have $\eps_{x,-\alpha}=e$. Therefore, we proved that $[\eps_x^{-1},y_{b_{-\beta}}^{-1}]$ has trivial component in $G^{(2)}_{\beta}(S)$. By Sublemma G we have $\eps_{x,0}=e$; hence the component of $[\eps_x^{-1},y_{b_{-\beta}}^{-1}]$ in $G^{(2)}_{\beta}(S)$ depends only on $\eps_{x,-\beta}^{-1}$. Therefore, $[\eps_{x,-\beta}^{-1},y_{b_{-\beta}}^{-1}]$ has trivial component in $G^{(2)}_{\beta}(S)$.
	
	Now let us prove that we have $(\eps_{x,-\beta})^{\tilde{w}_{-\beta}}\in G_{\beta}^{(2)}(S)$. First we prove that we have \mbox{$(\eps_{x,-\beta})^{\tilde{w}_{-\beta}}\in G_{\alpha}^{(2)}(S)G_{\beta}^{(2)}(S)G_{2\beta-\alpha}^{(2)}(S)$}. By \cite[Theorem~2, item~2]{VoronDiophantine} it it enough to show that $(\eps_{x,-\beta})^{\tilde{w}_{-\beta}}$ commutes with $G_{\gamma}^{(2)}(S)$ for all $\gamma$ such that $\gamma+\beta\notin\Phi\cup\{0\}$. It is the same that to show that $\eps_{x,-\beta}$ commutes with $\leftact{\tilde{w}_{-\beta}}{(G_{\gamma}^{(2)}(S))}$. By Sublemma E we have $\leftact{\tilde{w}_{-\beta}}{(G_{\gamma}^{(2)}(S))}\sub G_{s_{\beta}(\gamma)}^{(2)}(S)G(\tilde{\Phi},S,I)$. Since $\gamma+\beta\notin\Phi\cup\{0\}$, we have $s_{\beta}(\gamma)-\beta\notin\Phi\cup\{0\}$; hence $\eps_{x,-\beta}$ commutes with $G_{s_{\beta}(\gamma)}^{(2)}(S)$; and since $G(\tilde{\Phi},S,I)$ is abelian, $\eps_{x,-\beta}$ also commutes with $G(\tilde{\Phi},S,I)$. Therefore, we conclude that $(\eps_{x,-\beta})^{\tilde{w}_{-\beta}}\in G_{\alpha}^{(2)}(S)G_{\beta}^{(2)}(S)G_{2\beta-\alpha}^{(2)}(S)$. However, since $(\eps_{x,-\beta})^{\tilde{w}_{-\beta}}$ is a product of elements from $G_{\beta}^{(2)}(S)$ and $G_{-\beta}^{(2)}(S)$, it must have trivial components in $G_{\alpha}^{(2)}(S)$ and $G_{2\beta-\alpha}^{(2)}(S)$; hence, we have $(\eps_{x,-\beta})^{\tilde{w}_{-\beta}}\in G_{\beta}^{(2)}(S)$.
	
	Now it is easy to see that $(\eps_{x,-\beta})^{\tilde{w}_{-\beta}}$ must be equal to the component of $[\eps_{x,-\beta}^{-1},y_{b_{-\beta}}^{-1}]$ in $G_{\beta}^{(2)}(S)$, which we prove to be trivial. Therefore, we have $(\eps_{x,-\beta})^{\tilde{w}_{-\beta}}=e$; hence, $\eps_{x,-\beta}=e$. Using that and the fact that by Sublemma G we have $\eps_{x,-\alpha}=e$, we can conclude by Sublemma B that $\eps_{x,\alpha-2\beta}=e$.
	
	\underline{Subsublemma is proved.}

	\smallskip
	
	Now let us prove that $\eps_{b_{-\beta},\beta-\alpha}=e$. Consider two cases. The first case is where the map $\map{u}{\tilde{\Phi}\cup\{0\}}{\Phi\cup\{0\}}$ is not a bijection; hence we are assuming that $b_{-\beta}$ commutes with $b_{\beta-\alpha}$. Then Subsublemma implies that $\eps_{b_{\beta-\alpha},-\beta}=\eps_{b_{\beta-\alpha},\alpha-2\beta}=e$. Thus we prove that $\eps_{b_{\beta-\alpha},-\beta}=e$ based on the assumption that  $\eps_{b_{-\beta},2\beta-\alpha}=e$. Therefore, we can mirror the entire argument and use the fact that $\eps_{b_{\beta-\alpha},\alpha-2\beta}=e$ to prove that $\eps_{b_{-\beta},\beta-\alpha}=e$. The second case is  where $\map{u}{\tilde{\Phi}\cup\{0\}}{\Phi\cup\{0\}}$ is a bijection. Since $b_{-\beta}^{h_{\alpha}^{(1)}}=b_{-\beta}^{-1}$ we have
	$$
	y_{b_{-\beta}}^{-1}\eps_{b_{-\beta}}^{h_{\alpha}^{(2)}}=\theta(b_{-\beta})^{h_{\alpha}^{(2)}}=\theta(b_{-\beta}^{h_{\alpha}^{(1)}})=\theta(b_{-\beta}^{-1})=\theta(b_{-\beta})^{-1}=\eps_{b_{-\beta}}^{-1}y_{b_{-\beta}}^{-1}=y_{b_{-\beta}}^{-1}[y_{b_{-\beta}},\eps_{b_{-\beta}}^{-1}]\eps_{b_{-\beta}}^{-1}\tp
	$$
	Hence, we have $\eps_{b_{-\beta}}^{h_{\alpha}^{(2)}}=[y_{b_{-\beta}},\eps_{b_{-\beta}}^{-1}]\eps_{b_{-\beta}}^{-1}$. Now compare the components in $G^{(2)}_{2\beta-\alpha}(S)$ on the left hand side and on the right hand side of this equality. On the left hand side this component is trivial, and on the right hand side it only comes from $[y_{b_{-\beta}},\eps_{b_{-\beta},\beta-\alpha}^{-1}]$ (because $\eps_{b_{-\beta},-\alpha}=e$ by Sublemma B). Therefore, we proved that $[y_{b_{-\beta}},\eps_{b_{-\beta},\beta-\alpha}^{-1}]=e$. Since $u$ is a bijection it follows that all the root subgroups are 1-parametric and the commutator between $G^{(2)}_{\beta}$ and $G^{(2)}_{\beta-\alpha}$ correspond up to sign to the product of parameters multiplied by~$2$. By Sublemma E, the image of $y_{b_{-\beta}}$ in $G(\tilde{\Phi},S/I)$ is a component of a Weyl element; hence the parameter for $y_{b_{-\beta}}$ is invertible modulo $I$, and hence is invertible in $S$. Since~$2$ is also invertible in $S$, we then have $\eps_{b_{-\beta},\beta-\alpha}=e$.
	
	Thus in both cases we proved that $\eps_{b_{-\beta},\beta-\alpha}=e$. Hence, Subsublemma implies that for any $x\in G^{(1)}_{\alpha-\beta}$ we have $\eps_{x,-\beta}=\eps_{x,\alpha-2\beta}=e$. Since we proved it based on the on the assumption that  $\eps_{b_{-\beta},2\beta-\alpha}=e$, it follows that we can mirror the entire argument and use the newly established fact that $\eps_{b_{\beta-\alpha},\alpha-2\beta}=e$ to prove that for any $x\in G^{(1)}_{\beta}$ we have $\eps_{x,\beta-\alpha}=\eps_{x,2\beta-\alpha}=e$.
	
	Now we prove the statement of the Sublemma. Let $x\in G_{\alpha}^{(1)}$, we must prove that $\eps_{x,0}\in C^{\mathrm{us}}$. By \cite[Theorem~2, item~1]{VoronDiophantine} it is enough to prove that $\eps_{x,0}$ commutes with $G^{(2)}_{\gamma}(S)$ for all non-ultrashort $\gamma\in\Phi$ such that $\angle(\alpha,\gamma)\le\pi/2$. If $\gamma\perp\alpha$, then this follows from Sublemma~F. Now let us prove that $\eps_{x,0}$ commutes with $G^{(2)}_{\beta}(S)$. Let $x'\in G^{(1)}_{\beta}$. By Sublemma D we have $[y_x^{-1},\eps_{x'}^{-1}]=[y_{x'}^{-1},\eps_x^{-1}]$. Now compare the components in $G^{(2)}_{\beta}(S)$ on the left hand side and on the right hand side of this equality. On the left hand side this component is trivial because it only come from $[y_x^{-1},\eps_{x',\beta-\alpha}^{-1}]$ and we prove above that $\eps_{x',\beta-\alpha}=e$. On the right hand side this component only come from $[y_{x'}^{-1},\eps_{x,0}^{-1}]$ because by Sublemma~F we have $\eps_x=\eps_{x,0}\eps_{x,\alpha}$. Therefore, $\eps_{x,0}$ commutes with $y_{x'}$. Since all possible $y_{x'}$ generate $G^{(2)}_{\beta}(S)$ as an $S$-module, it follows that $\eps_{x,0}$ commutes with $G^{(2)}_{\beta}(S)$. Similarly we obtain that $\eps_{x,0}$ commutes with $G^{(2)}_{\alpha-\beta}(S)$. Further, it follows from Lemma~\ref{WeylForSplit} that $G^{(2)}_{\alpha}(S)\sub [G^{(2)}_{\alpha-\beta}(S),G^{(2)}_{2\beta-\alpha}(S)]G^{(2)}_{\beta}(S)$; hence, $\eps_{x,0}$ commutes with $G^{(2)}_{\alpha}(S)$. It remains to consider the case, where $\gamma$ is a short root such that $\angle(\alpha,\gamma)=\pi/4$ and distinct from $\beta$ and $\alpha-\beta$. In this case $\gamma-\beta$ is a short root orthogonal to $\alpha$. Thus it follows from Lemma~\ref{WeylForSplit} that $G^{(2)}_{\gamma}(S)\sub [G^{(2)}_{\beta}(S),G^{(2)}_{\gamma-\beta}(S)]$; hence, $\eps_{x,0}$ commutes with $G^{(2)}_{\gamma}(S)$.
	
	\underline{Sublemma H is proved.}
	
	\medskip
	
	Suppose that the root system $\Phi$ is of type $C$ or $BC$. For every long root $\alpha\in\Phi$ let $g_{-\alpha}\in G_{-\alpha}^{(2)}(S)\cap G(\tilde{\Phi},S,I)$ be as in Sublemma H. Let $g$ be the product of $g_{-\alpha}$ for all long roots $\alpha\in\Phi$. Now if we replace the homomorphism $\theta(-)$ with the homomorphism $\theta(-)^{g}$, then for every long root $\alpha$  we will have $\eps_{x,0}\in C^{\mathrm{us}}$ for every $x\in G_{\alpha}^{(1)}$, because factors of~$g$ other than $g_{-\alpha}$ do not affect $\eps_{x,0}$. Note that $\theta(h^{(1)}_{\gamma})$ has not changed for all $\gamma\in \Phi$; thus results of all the previous sublemmas still hold.
	
	In other words we may now assume without loss of generality that if $\Phi$ is of type $C$ or $BC$, the for every long root $\alpha\in\Phi$ and for every $x\in G_{\alpha}^{(1)}$ we have $\eps_{x,0}\in C^{\mathrm{us}}$. This assumption will be referred to as the \underline{correction assumption}.
	
	 \medskip
	
	\underline{Sublemma I.} Let $\alpha$,$\beta$,$\gamma\in\Phi$ be such that $\alpha$ is long, $\gamma\notin2\Phi$, $\beta$ is not ultrashort, $\alpha+\beta\notin\Phi\cup\{0\}$ and $\angle(\alpha,\gamma)>\pi/2$. Then we have $\eps_{x,\gamma}=e$ for every  $x\in G_{\beta}^{(1)}$.
	
	\medskip
	
	\underline{Proof of Sublemma I.}
	
	By \cite[Theorem~1]{VoronDiophantine} it is enough to prove that $\eps_x$ commutes with all the $y_{x'}$ for $x'\in G_{\alpha}^{(1)}$ (all such $y_{x'}$ generate $G_{\alpha}^{(1)}$ as an $S$-module). By Sublemma D we have $[y_x^{-1},\eps_{x'}^{-1}]=[y_{x'}^{-1},\eps_x^{-1}]$. It follows from Sublemma F and the correction assumption that $[y_x^{-1},\eps_{x'}^{-1}]=e$. Therefore, $[y_{x'}^{-1},\eps_x^{-1}]=e$, i.e.  $\eps_x$ commutes with $y_{x'}$ as required.
	
	\underline{Sublemma I is proved.}
	
	\medskip	
	
	\underline{Sublemma J.} Let $\alpha\in\Phi$ be a short root. Then we have $\theta(G_{\alpha}^{(1)})\sub G_{\alpha}^{(2)}(S)$.
	
	\medskip
	
	\underline{Proof of Sublemma J.}
	
	Let $x\in G_{\alpha}^{(1)}$. First let us prove that for every $\beta\in\Phi\sm (2\Phi)\cup\{\alpha\}$ we have $\eps_{x,\beta}=e$. If $\Phi$ is not of type $C$ or $BC$, then this follows from combination of Sublemmas A and I. If $\Phi$ is of type $C$ or $BC$, then by Sublemma I we have $\eps_{x,\beta}=e$ unless $\beta$ is a long or ultrashort root with $\angle(\alpha,\beta)=\pi/4$. Once this is established we can use combination of Sublemmas A and B to deal with these roots as well. Thus we have $\theta(x)=y_x\eps_{x,\alpha}\eps_{x,0}$.
	
	Now let $\gamma\in \Phi$ be a non-ultrashort root such that $\<\alpha,\gamma\>$ is odd. Therefore, we have
	$$
	\theta(x)^{-1}=\theta(x^{-1})=\theta(x^{h_{\gamma}^{(1)}})=\theta(x)^{h_{\gamma}^{(2)}}\tp
	$$
	Now using what we prove before, we can compare the components in $L^{(2)}(S)$ on the left hand side and on the right hand side of this equality and conclude that $\eps_{x,0}=\eps_{x,0}^{-1}$; hence, since $2\in S^*$ we have $\eps_{x,0}=e$.
	
	\underline{Sublemma J is proved.}
	
	\medskip	
	
	\underline{Sublemma K.} Let $\Phi$ be of type $BC$. Let $\alpha$ be a ultrashort root, and $x\in G_{\alpha}^{(1)}$. Then we have $\theta(x)=y_x\eps_{x,\alpha}\eps_{x,0}$, and $\eps_{x,0}\in C^{\mathrm{us}}$.
	
	\medskip
	
	\underline{Proof of Sublemma K.}
	
	By \cite[Theorem~2, item~1]{VoronDiophantine} it is enough to show that for every $\beta\in\Phi$ if $\alpha+\beta\notin \Phi\cup\{0\}$, then $\eps_x$ commutes with some subset of $G_{\beta}^{(2)}(S)$ that generates it as an $S$-module.
	
	First consider the case, where $\beta$ is short. Let us prove that $\eps_x$ commutes with all the $y_{x'}$, where $x'\in G_{\beta}^{(1)}$. By Sublemma D we have $[y_x^{-1},\eps_{x'}^{-1}]=[y_{x'}^{-1},\eps_x^{-1}]$.  It follows from Sublemma~J that $[y_x^{-1},\eps_{x'}^{-1}]=e$. Therefore, $[y_{x'}^{-1},\eps_x^{-1}]=e$, i.e.  $\eps_x$ commutes with $y_{x'}$ as required.
	
	Now consider the case, where $\beta$ is a long root such that $\alpha\perp\beta$. By Sublemma D we have $[y_x^{-1},\eps_{x'}^{-1}]=[y_{x'}^{-1},\eps_x^{-1}]$. It follows from Sublemma F that $[y_x^{-1},\eps_{x'}^{-1}]\in G_{\alpha}^{(2)}(S)$; however $[y_{x'}^{-1},\eps_x^{-1}]$ have trivial component in $G_{\alpha}^{(2)}(S)$ because $\alpha-\beta\notin\Phi\cup\{0\}$. Therefore, $[y_{x'}^{-1},\eps_x^{-1}]=e$ as required.
	
	The last case is where $\beta=2\alpha$. Take a short root $\beta'\in\Phi$ such that $\angle(\alpha,\beta)=\pi/4$. By previous cases $\eps_x$ commutes with $y_x'$ for all $x'\in G_{\beta'}^{(1)}$, for all $x'\in G_{2(\alpha-\beta')}^{(1)}$, and for all $x'\in G_{2\alpha-\beta'}^{(1)}$. Since $G(\tilde{\Phi},S,I)$ is abelian, we obtain that $\eps_x$ commutes with $\theta(G_{\beta'}^{(1)})$, with $\theta(G_{2(\alpha-\beta')}^{(1)})$ and with $\theta(G_{2\alpha-\beta'}^{(1)})$. Existence of the Weyl elements implies that any element of $G_{2\alpha}^{(1)}$ can be obtained as a commutator of element from $G_{\beta'}^{(1)}$ and $G_{2(\alpha-\beta')}^{(1)}$ taken modulo $G_{2\alpha-\beta'}^{(1)}$. Therefore, $\eps_x$ commutes with $\theta(G_{2\alpha}^{(1)})$; hence, $\eps_x$ commutes with $y_x'$ for all $x'\in G_{2\alpha}^{(1)}$.
	
	\underline{Sublemma K is proved.}
	
		\medskip	
	
	\underline{Sublemma L.} Let $\alpha\in\Phi$ be a long root. Then we have $\theta(G_{\alpha}^{(1)})\sub G_{\alpha}^{(2)}(S)$.
	
	\medskip
	
	\underline{Proof of Sublemma L.}
	
	By Sublemma F and the correction assumption, the statement is true unless $\Phi$ is of type $BC$. Now suppose that $\Phi$ is of type $BC$.
	
	Let $x\in G_{\alpha}^{(1)}$. By Lemma~\ref{EgorFormulas} it is enough to show that for every $\beta\in\Phi$ if $\alpha+\beta\notin \Phi\cup\{0\}$, then $\eps_x$ commutes with some subset of $G_{\beta}^{(2)}(S)$ that generates it as an $S$-module. If $\beta$ is not ultrashort, then this again follows from Sublemma F and the correction assumption. Now let $\beta$ be ultrashort. Let us prove that $\eps_x$ commutes with all the $y_{x'}$, where $x'\in G_{\beta}^{(1)}$. By Sublemma D we have $[y_x^{-1},\eps_{x'}^{-1}]=[y_{x'}^{-1},\eps_x^{-1}]$.  It follows from Sublemma K that we have $[y_x^{-1},\eps_{x'}^{-1}]=e$. Therefore, $[y_{x'}^{-1},\eps_x^{-1}]=e$, i.e.  $\eps_x$ commutes with $y_{x'}$ as required.
	
	\underline{Sublemma L is proved.}
	
	\medskip	
	
		\underline{Sublemma M.} Let $\alpha\in\Phi$ be a ultrashort root. Then we have $\theta(G_{\alpha}^{(1)})\sub G_{\alpha}^{(2)}(S)$.
	
	\medskip
	
	\underline{Proof of Sublemma M.}
	
	Let $x\in G_{\alpha}^{(1)}$. By Sublemma K it is enough to show that $\eps_{x,0}=e$. Let $\gamma\in \Phi$ be a non-ultrashort root such that $\<\alpha,\gamma\>$ is odd. Then we have $x^{h_{\gamma}^{(1)}}x\in G_{2\alpha}^{(1)}$. Hence, by Sublemma~L we have $\theta(x^{h_{\gamma}^{(1)}}x)\in G_{2\alpha}^{(2)}(S)$. On the other hand we have
	$$
	\theta(x^{h_{\gamma}^{(1)}}x)=\theta(x)^{h_{\gamma}^{(2)}}\theta(x)\in \eps_{x,0}^2 G_{2\alpha}^{(2)}(S)\tp
	$$
	Therefore, we have $\eps_{x,0}^2=e$; hence, since $2\in S^*$ we have $\eps_{x,0}=e$.
	
	\underline{Sublemma M is proved.}
	
	\medskip	
	
	The Lemma now follows from Sublemmas J, L and M.
\end{proof}

{\prop\label{InfinitesiamiPartI} Under the umbrella assumptions of this Section, there exist an element $g\in G(\tilde{\Phi},S,I)$ such that for the homomorphism $\theta'(x)=\theta(x)^g$ we have $\theta'(h^{(1)}_{\alpha})=h^{(2)}_\alpha$ for all non-ultrashort $\alpha\in\Phi$; and  $\theta'(G^{(1)}_{\alpha})\sub G^{(2)}_{\alpha}(S)$ for all $\alpha\in\Phi$.}

\begin{proof}
	By assumption we have $\theta(h_{\alpha}^{(1)})=h_{\alpha}^{(2)}\eps_{\alpha}$, where $\eps_\alpha\in G(\tilde{\Phi},S,I)$.  Further we can decompose the element $\eps_{\alpha}$ as $\eps_{\alpha}=\prod_{\gamma\in\Phi\sqcup\{0\}\sm 2\Phi} \eps_{\alpha,\gamma}$, where $\eps_{\alpha,0}\in L^{(2)}(S)\cap G(\tilde{\Phi},S,I)$ and $\eps_{\alpha,\gamma}\in G^{(2)}_{\gamma}(S)\cap G(\tilde{\Phi},S,I)$ for all $\gamma\in \Phi\sm 2\Phi$. Let us prove that there exist an element $g'\in G(\tilde{\Phi},S,I)$ such that after we replace the homomorphism $\theta(-)$ with the homomorphism $\theta(-)^{g'}$ we will have $\eps_{\alpha,\gamma}=e$ for all non-ultrashort $\alpha\in\Phi$ and all non-ultrashort $\gamma\in\Phi\sm2\Phi$, and $\eps_{\alpha,\gamma}\in G^{(2)}_{2\gamma}(S)$ for all non-ultrashort $\alpha\in\Phi$ and all ultrashort $\gamma\in\Phi$.
	
	We will look for $g'$ in a form $g'=\prod_{\beta\in\Phi\sm2\Phi} g'_{\beta}$, where 	$g'_{\beta}\in G^{(2)}_{\beta}(S)\cap G(\tilde{\Phi},S,I)$. The requirement on $g'$ is that for every $\gamma\in\Phi\sm 2\Phi$ and for every non-ultrashort $\gamma\in\Phi$ the component of $h_{\alpha}^{(2)}(h_{\alpha}^{(2)}\eps_{\alpha})^{g'}$ in $G^{(2)}_{\gamma}(S)$ must be either trivial, or belong to $G^{(2)}_{2\gamma}(S)$. When we rewrite this in terms of $\eps_{\alpha,\gamma}$ and $g'_{\beta}$ we get two type of equation. The first type is that whenever $\<\gamma,\alpha\>$ is even we must have $\eps_{\alpha,\gamma}=e$. The second type is that whenever $\<\gamma,\alpha\>$ is odd we must have $(g'_{\gamma})^2=\eps_{\alpha,\gamma}^{-1}$ or $(g'_{\gamma})^2\in \eps_{\alpha,\gamma}^{-1}G^{(2)}_{2\gamma}(S)$ if $\gamma$ is ultrashort.
	
	Let us show that the equations of the first type hold true automatically. Indeed, since $(h_{\alpha}^{(1)})^2=e$, it follows that $(h_{\alpha}^{(2)}\eps_{\alpha})^2=e$. Looking at the component in $G^{(2)}_{\gamma}(S)$, when $\<\gamma,\alpha\>$ is even, we get $\eps_{\alpha,\gamma}^2=e$; and since $2\in S^*$, it follows that $\eps_{\alpha,\gamma}=e$.
	
	For the equations of the second type we must prove that there exists $g'$ that satisfies all of them simultaneously. Since $2\in S^*$ we have no problem with extracting square roots of $\eps_{\alpha,\gamma}^{-1}$ (modulo $G^{(2)}_{2\gamma}(S)$ if $\gamma$ is ultrashort). Thus we only have to prove that for any non-ultrashort $\alpha_1$,$\alpha_2\in\Phi$ and any $\gamma\in\Phi\sm2\Phi$ such that $\<\gamma,\alpha_1\>$ and $\<\gamma,\alpha\>$ are odd we have $\eps_{\alpha_1,\gamma}=\eps_{\alpha_2,\gamma}$ (or $\eps_{\alpha_1,\gamma}\in \eps_{\alpha_2,\gamma}G^{(2)}_{2\gamma}(S)$), so that the equations do not contradict each other.
	
	Since we have $h_{\alpha_1}^{(1)}h_{\alpha_2}^{(1)}=h_{\alpha_2}^{(1)}h_{\alpha_1}^{(1)}$, it follows that $h_{\alpha_1}^{(2)}\eps_{\alpha_1}h_{\alpha_2}^{(2)}\eps_{\alpha_2}=h_{\alpha_2}^{(2)}\eps_{\alpha_2}h_{\alpha_1}^{(2)}\eps_{\alpha_1}$, i.e $\eps_{\alpha_1}^{h_{\alpha_2}^{(2)}}\eps_{\alpha_2}=\eps_{\alpha_2}^{h_{\alpha_1}^{(2)}}\eps_{\alpha_1}$. Comparing component in $G^{(2)}_{\gamma}(S)$, we obtain $\eps_{\alpha_1,\gamma}^2=\eps_{\alpha_2,\gamma}^2$ (or $\eps_{\alpha_1,\gamma}^2\in \eps_{\alpha_2,\gamma}^2G^{(2)}_{2\gamma}(S)$); and since $2\in S^*$, it follows that $\eps_{\alpha_1,\gamma}=\eps_{\alpha_2,\gamma}$ (or $\eps_{\alpha_1,\gamma}\in \eps_{\alpha_2,\gamma}G^{(2)}_{2\gamma}(S)$) as required.
	
	Therefore, in the view of the above, we may assume that $\eps_{\alpha,\gamma}=e$ for all non-ultrashort $\alpha\in\Phi$ and all non-ultrashort $\gamma\in\Phi\sm2\Phi$, and $\eps_{\alpha,\gamma}\in G^{(2)}_{2\gamma}(S)$ for all non-ultrashort $\alpha\in\Phi$ and all ultrashort $\gamma\in\Phi$. Tis implies that $\eps_{\alpha}$ commutes with $h_{\alpha}^{(2)}$. Hence, the equation $(h_{\alpha}^{(2)}\eps_{\alpha})^2=e$ implies that $\eps_{\alpha}^2=e$; hence, $\eps_{\alpha}=e$; i.e $\theta(h^{(1)}_{\alpha})=h^{(2)}_\alpha$. The proposition now follows from Lemma~\ref{GiantLemma}.
	\end{proof}

\section{The scheme of adjustments}
\label{AdjustmentSec}

In this section we are assuming the assumptions of Theorem~\ref{main}; and additionally we are assuming that $\tilde{\Phi_1}=\tilde{\Phi_2}=\tilde{\Phi}$.

\medskip

Consider the following functor from the category of $R_2$-algebras into the category of sets:
$$
\Pin(S)=\{f\colon (G_2)_S \toiso G(\tilde{\Phi},-)_S\mid f \text{ is an isomorphism of group schemes over }S\}\tc
$$
where $(G_2)_S$ is the base change of $G_2$ to $S$, and $G(\tilde{\Phi},-)_S$ is the adjoint Chevalley--Demazure scheme over $S$.

{\lem\label{PinGood} The functor $\Pin(-)$ is smooth affine scheme of finite type over $R_2$.}

\begin{proof}
	By \cite[Exp. XXIV, Section 7]{SGAIII} there are affine schemes that classify group-scheme homomorphisms from $G_2$ to $G(\tilde{\Phi},-)$ and from $G(\tilde{\Phi},-)$ to $G_2$. The functor $\Pin(-)$ can be realized as a subscheme in their product defined by a closed condition that the product of homomorphisms in both orders is the identity homomorphism. Therefore, $\Pin(-)$ is an affine scheme.
	
	Now, clearly, $\Pin(-)$ is a torsor of a group scheme that classify group-scheme automorphisms of $G(\tilde{\Phi},-)$, which as a scheme is isomorphic to the product of $G(\tilde{\Phi},-)$ with a finite constant scheme. Thus $\Pin(-)$ is a torsor of a smooth group scheme of finite type, and hence is smooth and of finite type.
	\end{proof}

Now let $(T^{(1)},\Phi)$ be an isotropic pinning on $G_1$ that satisfies the assumption (a)--(d) in the second bullet in the statement of Theorem~\ref{main}; and let $\map{u}{\tilde{\Phi}\cup\{0\}}{\Phi\cup\{0\}}$ be the corresponding map. Let $G_{\alpha}^{(1)}=G_{\alpha}^{(1)}(R_1)$ for $\alpha\in\Phi$ be the corresponding root subgroups; and $h_{\alpha}^{(1)}$ are the elements from Lemma~\ref{haExist}. Let $(T^{(2)},\Phi)$ be an isotropic pre-pinning on $G(\tilde{\Phi},-)$ that is contained in the pinning and such that the corresponding map ${\tilde{\Phi}\sqcup\{0\}}\to{\Phi\sqcup\{0\}}$ is $u$. Note that by Lemma~\ref{WeylForSplit} this isotropic pre-pinning is actually an isotropic pinning. Let $G_{\alpha}^{(2)}$ for $\alpha\in\Phi$ be the corresponding root subgroup-subscheme (i.e. $G_{\alpha}^{(2)}=\prod_{u(\beta)\in\{\alpha,2\alpha\}} U_{\beta}$), and $h_{\alpha}^{(2)}$ be the elements from Lemma~\ref{haExist}.

For every $\alpha\in\Phi$ choose a finite set $X_{\alpha}\sub G_{\alpha}^{(1)}(R_1)$ such that $X_{\alpha}$ generate $G_{\alpha}^{(1)}(R_1)$ or $G^{(1)}_{\alpha}(R_1)/G^{(1)}_{2\alpha}(R_1)$ if $\alpha$ is ultrashort) as an $R_1$-module.

{\defn Let $S$ be an $R_2$-algebra. We call an element $f\in\Pin(S)$ an \underline{adjustment} if, firstly, we have $f_S(\theta(h_{\alpha}^{(1)}))=h_{\alpha}^{(2)}$ for all non-ultrashort $\alpha\in\Phi$; secondly, we have $f_S(\theta(G_{\alpha}^{(1)}))\sub G_{\alpha}^{(2)}(S)$ for all $\alpha\in\Phi$; and thirdly, $f_S(\theta(X_{\alpha}))$ generate $G_{\alpha}^{(2)}(S)$ (or $G^{(2)}_{\alpha}(S)/G^{(2)}_{2\alpha}(S)$ if $\alpha$ is ultrashort) as an $S$-module for all $\alpha\in\Phi$.}

{\lem\label{AdjustXIsEnough} Let $S$ be an $R_2$-algebra. Let $f\in\Pin(S)$ be such that firstly, $f_S(\theta(h_{\alpha}^{(1)}))=h_{\alpha}^{(2)}$ for all non-ultrashort $\alpha\in\Phi$; secondly, $f_S(\theta(X_{\alpha}))\sub G_{\alpha}^{(2)}(S)$ for all $\alpha\in\Phi$; and thirdly, $f_S(\theta(X_{\alpha}))$ generate $G_{\alpha}^{(2)}(S)$ {\rm (}or $G^{(2)}_{\alpha}(S)/G^{(2)}_{2\alpha}(S)$ if $\alpha$ is ultrashort{\rm )} as an $S$-module for all $\alpha\in\Phi$. Then $f$ is an adjustment.}
\begin{proof}
	It is enough to how that $f_S(\theta(G_{\alpha}^{(1)}))\sub G_{\alpha}^{(2)}(S)$. We can do it by showing that elements from $f_S(\theta(G_{\alpha}^{(1)}))$ satisfy the formulas form Lemma~\ref{EgorFormulas} with $f_S(\theta(X_{\gamma}))$ as generating set for $G_{\gamma}^{(2)}(S)$ for every $\gamma\in\Phi$. These formulas easily follow from the fact that $f_S\circ\theta$ is a group homomorphism.
	\end{proof}

\smallskip

Consider the following functor from the category of $R_2$-algebras into the category of sets:
$$
\Adjust(S)=\{f\in\Pin(S)\mid f\text{ is an adjustment}\}\tp
$$

{\lem\label{AdjustIsOfFiniteType} The functor $\Adjust(-)$ is a locally finitely presented quasi-affine scheme of finite type over $R_2$.}

\begin{proof}
	By Lemma~\ref{PinGood} the scheme $\Pin(-)$ is smooth, so in particular it is locally finitely presented. The functor $\Adjust(-)$ is a subfunctor of $\Pin(-)$ defined by conditions in Lemma~\ref{AdjustXIsEnough}. The conditions $f_S(\theta(h_{\alpha}^{(1)}))=h_{\alpha}^{(2)}$ and $f_S(\theta(X_{\alpha}))\sub G_{\alpha}^{(2)}(S)$ are clearly equivalent to certain finite collection of regular functions on $\Pin(-)$ vanishing in $f$; hence they define a locally finitely presented closed subscheme of finite type. Further for any $\alpha\in\Phi$ we can write a matrix with regular functions on $\Pin(-)$ as entries, such that the condition ''$f_S(\theta(X_{\alpha}))$ generate $G_{\alpha}^{(2)}(S)$ (or $G^{(2)}_{\alpha}(S)/G^{(2)}_{2\alpha}(S)$)'' is equivalent to maximal order minors of that matrix, when evaluates in $f$, generate a unit ideal of $S$. Therefore, this condition define an open subscheme, whose complement (in sense of schemes) is an intersection of closed subschemes each defined by a condition ''particular minor of the matrix vanishes''. The functor $\Adjust(-)$ is then the intersection of all these open subschemes; so, clearly it is a locally finitely presented quasi-affine scheme of finite type.
\end{proof}

{\lem\label{AdjustOverFields} For any maximal ideal $J\unlhd R_2$ there is a finite field extension $R_2/J \le K$ such that $\Adjust(K)\ne\emp$.}

\begin{proof}
	Let $I\unlhd R_1$ be a maximal ideal that correspond to $I$ by correspondence from Lemma~\ref{MaxCorr}. Then $\theta$ induces an isomorphism $\overline{\theta}\colon E_1(R_1/I)\toiso E_2(R_2/J)$. By Proposition~\ref{BT} $\overline{\theta}$ comes from a field isomorphism $\ph\colon R_1/I\toiso R_2/J$ and an $R_2/J$-group-scheme isomorphism $\overline{\Theta}\colon \phan^\ph (G_1)_{R_1/I}\toiso (G_2)_{R_2/J}$ (all the exceptional cases are excluded by assumptions of Theorem~\ref{main} and the assumption $\tilde{\Phi_1}=\tilde{\Phi_2}$). In particular, there is an isotropic pinning on the scheme $(G_2)_{R_2/J}$, with the same map $\map{u}{\tilde{\Phi}\cup\{0\}}{\Phi\cup\{0\}}$ as above, and such that  $\overline{\theta}$ maps the image of the root subgroups $G_{\alpha}^{(1)}(R_1)$ onto the corresponding root subgroups. The elements $h_{\alpha}^{(1)}$ must be mapped to the corresponding elements $h_{\alpha}$. Now take $K$ to be such extension that after a scalar extension to $K$ our isotropic pinning becomes contained in a pinning. Then it is easy to see that we have  $\Adjust(K)\ne\emp$.
\end{proof}

{\lem\label{AdjustIsSmooth} The scheme $\Adjust(-)$ is smooth over $R_2$.}
\begin{proof}
	By Lemma~\ref{AdjustIsOfFiniteType} the scheme $\Adjust(-)$ is locally finitely presented. Hence, by \cite[Ch I, \S 4, Item 4.5]{DemazureGabriel} it is enough to show that for any $R_2$ algebra $S$ and any ideal $I\unlhd S$ such that $I^2=0$ the reduction map $\Adjust(S)\to\Adjust(S/I)$ is surjective.
	
	Now let $\overline{f}\in \Adjust(S/I)$. Since by Lemma~\ref{PinGood} the scheme $\Pin(-)$ is smooth, there is $f'\in \Pin(S)$ such that its reduction modulo $I$ is $\overline{f}$. Using Proposition~\ref{InfinitesiamiPartI} for the map $f'_S\circ\theta$ (here all the umbrella assumptions of Section~\ref{InfinitesimalSecI} are satisfied because $\overline{f}$ is an adjustment) we conclude that there exists $g\in G(\tilde{\Phi},S,I)$ such that $f=(f')^g$ is an adjustment. Clearly, conjugation by $g$ does not affect reduction modulo $I$; hence the reduction of $f$ is $\overline{f}$.
\end{proof}

{\prop\label{AdjustExists} There is a faithfully flat finitely presented $R_2$-algebra $S$ such that $\Adjust(S)\ne\emp$.}
\begin{proof}
Let $\{\Spec S_i\}_{i=1}^n$ be an open cover of $\Adjust(-)$ by affine subschemes (here we used Lemma~\ref{AdjustIsOfFiniteType}). Take $S=\prod_{i=1}^n S_i$. The composition of morphisms
$$
\Spec S\to \Adjust(-)\to \Spec R_2
$$
defines a structure of $R_2$ algebra over $S$ in such a way that $\Adjust(S)\ne\emp$ by construction. It also follows from Lemma~\ref{AdjustIsOfFiniteType} that $S$ is finitely presented over $R_2$. The morphism $\Spec S\to \Adjust(-)$ by construction is flat and surjective on the underlying topological spaces. The morphism $\Adjust(-)\to \Spec R_2$ if flat by Lemma~\ref{AdjustIsSmooth}; hence, the map on the underlying topological spaces has an open image. By Lemma~\ref{AdjustOverFields} this image contains all the closed points. Since $\Spec R_2$ is affine this implies that the map on the underlying topological spaces is surjective. Therefore, the morphism $\Spec S\to \Spec R_2$ is flat and surjective on the underlying topological spaces. Hence $S$ is faithfully flat over~$R_2$.
\end{proof}

\section{Isotropic pinning is preserved}
\label{PinningSec}

In this section we are assuming the assumptions of Theorem~\ref{main}; and additionally we are assuming that $\tilde{\Phi_1}=\tilde{\Phi_2}=\tilde{\Phi}$.

Let $(T^{(1)},\Phi)$ be an isotropic pinning on $G_1$ that satisfies the assumption (a)--(d) in the second bullet in the statement of Theorem~\ref{main}; and let $\map{u}{\tilde{\Phi}\cup\{0\}}{\Phi\cup\{0\}}$ be the corresponding map. Let $G_{\alpha}^{(1)}=G_{\alpha}^{(1)}(R_1)$ for $\alpha\in\Phi$ be the corresponding root subgroups; and $h_{\alpha}^{(1)}$ be the elements from Lemma~\ref{haExist}. Let $(T^{(\spl)},\Phi)$ be an isotropic pinning on $G(\tilde{\Phi},-)$ that is contained in the pinning and such that the corresponding map ${\tilde{\Phi}\sqcup\{0\}}\to{\Phi\sqcup\{0\}}$ is $u$. Let $G_{\alpha}^{(\spl)}$ for $\alpha\in\Phi$ be the corresponding root subgroup-subscheme and $h_{\alpha}^{(\spl)}$ be the elements from Lemma~\ref{haExist}.

For every root $\alpha\in\Phi$ define a subfunctor $G_{\alpha}^{(2)}(-)$ of the functor $G_2(-)$ in the following way: for any $R_2$-algebra $A$
$$
G_{\alpha}^{(2)}(A)=\{g\in G_2(A)\colon \text{(F'1)--(F'3) hold}\}\tc
$$
where

(F'1) $[g,\theta(x_{\gamma})]=e$ for every $\gamma\in\Phi$ such that $\alpha+\gamma\notin\Phi\cup\{0\}$ and every $x_{\gamma}\in G_{\gamma}^{(1)}$.

(F'2) If $\alpha$ is ultrashort, then $[[g,\theta(x_{\alpha})],\theta(x_{\gamma})]=e$ for every root $\gamma\in\Phi$ such that $2\alpha+\gamma\notin\Phi\cup\{0\}$ and every $x_{\gamma}\in G_{\gamma}^{(1)}$ and $x_{\alpha}\in G_{\alpha}^{(1)}$.

(F'3) If $\alpha$ is a short root and $\Phi$ is of type $C$ or $BC$ (including $C_2=B_2$), and $\alpha'$ and $\alpha''$ are distinct long roots such that $\angle(\alpha,\alpha')=\angle(\alpha,\alpha'')=\pi/4$, then $[[g,\theta(x_{\alpha'-\alpha})],\theta(x_{\gamma})]=e$ for every $\gamma\in\Phi$ such that $\alpha'+\gamma\notin\Phi\cup\{0\}$ and every $x_{\gamma}\in G_{\gamma}^{(1)}$ and $x_{\alpha'-\alpha}\in G_{\alpha'-\alpha}^{(1)}$; and $[[g,\theta(x_{\alpha''-\alpha})],\theta(x_{\gamma})]=e$ for every $\gamma\in\Phi$ such that $\alpha''+\gamma\notin\Phi\cup\{0\}$ and every $x_{\gamma}\in G_{\gamma}^{(1)}$ and $x_{\alpha''-\alpha}\in G_{\alpha''-\alpha}^{(1)}$.

\smallskip

Here we identified $\theta(x)$ for $x\in E_1(R_1)$ with the image of $\theta(x)$ in $G_2(A)$ under the group homomorphism induced by the structural homomorphism $R_2\to A$.

{\prop\label{PinningPreserved} There is an isotropic pinning $(T^{(2)},\Phi)$ on $G_2$ with $G_{\alpha}^{(2)}(-)$ as root subgroup-subschemes; and we have $\theta(G_{\alpha}^{(1)}(R_1))=G_{\alpha}^{(2)}(R_2)$.}

\begin{proof}
	It can be easily seen from the definition that $G_{\alpha}^{(2)}(-)$ are the closed subschemes of $G_2$. By Proposition~\ref{AdjustExists} there is a faithfully flat finitely presented $R_2$-algebra $S$ and a group-scheme isomorphism $f\colon (G_2)_S \toiso G(\tilde{\Phi},-)_S$ such that $f$ is an adjustment. It follows from Lemma~\ref{EgorFormulas} applied to $G(\tilde{\Phi},-)_S$ that $f$ turns $(G_{\alpha}^{(2)})_S$ into $G_{\alpha}^{(\spl)}$. For a closed subscheme, the property of being a subgroup-subscheme is fppf-local; hence, $G_{\alpha}^{(2)}$ are subgroup-subschemes. For any non-ultrashort root $\alpha\in\Phi$ let $\g_\alpha=\Lie(G_{\alpha}^{(2)})$ be the Lie algebra of $G_{\alpha}^{(2)}$. If $\alpha$ is ultrashort, then denote by $\g_{\alpha}$ the submodule of all $\Lie(G_{\alpha}^{(2)})$ defined as $\{x\in \Lie(G_{\alpha}^{(2)})\colon [x,\g_{\beta}]\sub \Lie(G_{\alpha+\beta}^{(2)})\}$, where $\beta$ is a short root such that $\angle(\alpha,\beta)=3\pi/4$. It is easy to see that $f$ turns $g_{\alpha}\otimes S$ into the root subspaces of the isotropic pinning $(T^{(\spl)},\Phi)$. Hence, by flat descent we can conclude that $\g_{\alpha}$ are the direct summands of $\Lie(G_2)$ and as a projective modules have constant rank.
	
	Now define the subfunctor $T^{(2)}(-)$ of the functor $G_2(-)$ in the following way: for any $R_2$-algebra $A$ define $T^{(2)}(A)$ as a subset of $G_2(A)$ that consists of elements $g\in G_2(A)$ such that the adjoint action of $g$ preserves all the $\g_{\alpha}\otimes A$ and acts on each such algebra by multiplication by a scalar (the scalar may depend on $\alpha$). It follows from what we showed above that $T^{(2)}$ is a closed subscheme of $G_2$ and $f$ turns $T^{(2)}$ into $T^{(\spl)}$. Therefore, $T^{(2)}$ is a torus; and if $\Delta\sub\Phi$ is a basis, then scalars that $T^{(2)}$ acts on $\g_{\alpha}$ for $\alpha\in\Delta$ with define an isomorphism $T^{(2)}\simeq \mathbb{G}_m^{\rk\Phi}$. Now, clearly, $(T^{(2)},\Phi)$ is an isotropic pre-pinning on $G_2$ with $G_{\alpha}^{(2)}$ as root subgroup-subschemes.
	
	The fact that $\theta(G_{\alpha}^{(1)}(R_1))=G_{\alpha}^{(2)}(R_2)$ follows from Lemma~\ref{EgorFormulas} applied to $G_1$ and the fact that $\theta$ is an isomorphism.
	
	Now it remains to show that $(T^{(2)},\Phi)$ is an isotropic pinning. Let $w_{\alpha}^{(1)}\in G_{\alpha}^{(1)}G_{-\alpha}^{(1)}G_{\alpha}^{(1)}$ be the Weyl elements for the pinning $(T^{(1)},\Phi)$ on $G_1$. We will prove that $\theta(w_{\alpha}^{(1)})$ are the Weyl elements for $(T^{(2)},\Phi)$. First it follows from what we proved previously that we have $\theta(w_{\alpha}^{(1)})\in G_{\alpha}^{(2)}(R_2)G_{-\alpha}^{(2)}(R_2)G_{\alpha}^{(2)}(R_2)$. Since $\theta$ is an isomorphism, it follows that $G_{\beta}^{(2)}(R_2)^{\theta(w_{\alpha}^{(1)})}=G_{s_{\alpha}(\beta)}^{(2)}(R_2)$. Now if $A$ is an $R_2$-algebra, then, clearly, the image of $G_{\beta}^{(2)}(R_2)$ generate $G_{\beta}^{(2)}(A)$ (or $G_{\beta}^{(2)}(A)/G_{2\beta}^{(2)}(A)$ if $\beta$ is ultrashort) as an $A$-module. Since by Lemma~\ref{EgorFormulas} the subgroups $G_{\beta}^{(2)}(A)$ can be defined through such generating subsets, it follows that we have $G_{\beta}^{(2)}(A)^{\theta(w_{\alpha}^{(1)})}=G_{s_{\alpha}(\beta)}^{(2)}(A)$. Therefore, $\theta(w_{\alpha}^{(1)})$ are the Weyl elements.
\end{proof}

\section{Isomorphism of rings}
\label{RingsSec}

In this section we are assuming the assumptions of Theorem~\ref{main}; and additionally we are assuming that $\tilde{\Phi_1}=\tilde{\Phi_2}=\tilde{\Phi}$. Let $(T^{(1)},\Phi)$ be an isotropic pinning on $G_1$ that satisfies the assumption (a)--(d) with the corresponding root subgroups $G_{\alpha}^{(1)}(R_1)$. Let $(T^{(2)},\Phi)$ be an isotropic pinning on $G_2$ from the Proposition~\ref{PinningPreserved} with the corresponding root subgroups $G_{\alpha}^{(2)}(R_2)$. Thus we have $\theta(G_{\alpha}^{(1)}(R_1))=G_{\alpha}^{(2)}(R_2)$.

{\prop\label{RingIsom} There is a ring isomorphism $\ph\colon R_1\toiso R_2$ such that for any $\alpha\in\Phi$, for any $x\in G_{\alpha}^{(1)}(R_1)$ and for any $\xi\in R_1$ we have $\theta(\xi\cdot x)=\ph(\xi)\cdot \theta(x)$ if $\alpha$ is non-ultrashort; and $\theta(\xi\cdot x)\in(\ph(\xi)\cdot \theta(x))G_{2\alpha}^{(2)}(R_2)$ if $\alpha$ is ultrashort.}
\begin{proof}
The proposition follows immediately from the proof of \cite[Theorem~5]{VoronDiophantine} and the fact that $\theta$ is an isomorphism.
\end{proof}

\section{Coincidence on ultrashort root elements}

{\prop \label{UltrashortFix} Let $G_1$ and $G_2$ be absolutely simple adjoint group schemes over rings $R_1$ and $R_2$, each with a common root datum of the geometric fibers. Let $\tilde{\Phi}$ be the absolute root system for both $G_1$ and $G_2$. Let $(T^{(1)},\Phi)$ and $(T^{(2)},\Phi)$ be isotropic pinnings on $G_1$ and $G_2$. Let $\map{u}{\tilde{\Phi_1}\cup\{0\}}{\Phi\cup\{0\}}$ be the map corresponding to both pinnings; and suppose that this map comes from a Tits index and that $\rk\Phi\ge 2$. Let $G_{\alpha}^{(1)}$ for $\alpha\in\Phi$ be the corresponding root subgroup-subschemes of $G_1$; and $G_{\alpha}^{(2)}$ for $\alpha\in\Phi$ be the corresponding root subgroup-subschemes of $G_2$. Let $S$ be an $R_2$ algebra. Let $\map{\theta,\theta'}{E_1(R_1)}{G_2(S)}$ be group homomorphisms such that
	
	$\bullet$ $\theta(G_{\alpha}^{(1)}(R_1))\sub G_{\alpha}^{(2)}(S)$ for every $\alpha\in\Phi$
	
	$\bullet$ $\theta(G_{\alpha}^{(1)}(R_1))$ generate $G_{\alpha}^{(2)}(S)$ (or  $G_{\alpha}^{(2)}(S)/G_{2\alpha}^{(2)}(S)$ if $\alpha$ is ultrashort) as an $S$-module.
	
	$\bullet$ $\restr{\theta}{G_{\alpha}^{(1)}(R_1)}=\restr{\theta'}{G_{\alpha}^{(1)}(R_1)}$ for all non-ultrashort $\alpha\in\Phi$
	
	$\bullet$ for any ultrashort $\alpha\in\Phi$ and any $x\in G_{\alpha}^{(1)}(R_1)$ we have $\theta(x)\in\theta'(x)G_{2\alpha}^{(2)}(S)$.
	
	Then we have $\theta=\theta'$.
}
\begin{proof}
	Clearly, the only case, where we have something to prove, is where $\Phi$ has type $BC$; and we only have to prove that $\restr{\theta}{G_{\alpha}^{(1)}(R_1)}=\restr{\theta'}{G_{\alpha}^{(1)}(R_1)}$ for all ultrashort $\alpha\in\Phi$.
	
	Let $\alpha\in\Phi$ be an ultrashort root; and let $x\in G_{\alpha}^{(1)}(R_1)$. Let $\beta$ be a short root such that $\angle(\alpha,\beta)=3\pi/4$. Then for any $y\in G_{\beta}^{(1)}(R_1)$ we have
	\begin{multline*}
	[\theta(x)\theta'(x)^{-1},\theta(y)]=[[\theta'(x)^{-1},\theta(y)],\theta(x)]\cdot[\theta'(x)^{-1},\theta(y)]\cdot [\theta(x),\theta(y)]=[[\theta'(x)^{-1},\theta'(y)],\theta(x)]\cdot\\ \cdot[\theta'(x)^{-1},\theta'(y)]\cdot [\theta(x),\theta(y)]=[\theta'([x^{-1},y]),\theta(x)] \cdot\theta'([x^{-1},y])\cdot [\theta(x),\theta(y)]\in G_{2(\alpha+\beta)}^{(2)}(S)\cdot\\ \cdot[\theta([x^{-1},y]),\theta(x)] \cdot\theta([x^{-1},y])\cdot [\theta(x),\theta(y)]=G_{2(\alpha+\beta)}^{(2)}(S)\theta([[x^{-1},y],x] \cdot[x^{-1},y]\cdot [x,y])=\\=G_{2(\alpha+\beta)}^{(2)}(S)\theta([xx^{-1},y])=G_{2(\alpha+\beta)}^{(2)}(S)\tp
	\end{multline*}
	
	Now since we have $\theta(x)\theta'(x)^{-1}\in G_{2\alpha}^{(2)}(S)$, and since $\theta(G_{\beta}^{(1)}(R_1))$ generate $G_{\beta}^{(2)}(S)$ as an $S$-module, it follows by \cite[Lemma~3]{VoronDiophantine} that $\theta(x)\theta'(x)^{-1}=e$.
\end{proof}

\section{The scheme of isomorphisms}
\label{IsomSec}

In this section we are assuming the assumptions of Theorem~\ref{main}; and additionally we are assuming that $\tilde{\Phi_1}=\tilde{\Phi_2}=\tilde{\Phi}$. Let $(T^{(1)},\Phi)$ be an isotropic pinning on $G_1$ that satisfies the assumption (a)--(d) with the corresponding root subgroups $G_{\alpha}^{(1)}(R_1)$. Let $(T^{(2)},\Phi)$ be an isotropic pinning on $G_2$ from the Proposition~\ref{PinningPreserved} with the corresponding root subgroups-subschemes $G_{\alpha}^{(2)}$. Thus we have $\theta(G_{\alpha}^{(1)}(R_1))=G_{\alpha}^{(2)}(R_2)$. Further let $\ph\colon R_1\toiso R_2$ be a ring isomorphism from Proposition~\ref{RingIsom}.

\smallskip

Consider the following functor from the category of $R_2$-algebras into the category of sets:
$$
\Isom(S)=\{\Theta\colon (\leftact{\ph}{G_1})_S\toiso (G_2)_S\mid \Theta \text{ is an isomorphism of group schemes over }S\}\tc
$$

where $(\leftact{\ph}{G_1})_S$ and $(G_2)_S$ are the base changes of the corresponding group-schemes to $S$.

{\lem\label{IsomGood} The functor $\Isom(-)$ is smooth affine scheme of finite type over $R_2$.}

\begin{proof}
	By \cite[Exp. XXIV, Section 7]{SGAIII} there are affine schemes that classify group-scheme homomorphisms from $\leftact{\ph}{G_1}$ to $G_2$ and from $G_2$ to $\leftact{\ph}{G_1}$. The functor $\Isom(-)$ can be realized as a subscheme in their product defined by a closed condition that the product of homomorphisms in both orders is the identity homomorphism. Therefore, $\Isom(-)$ is an affine scheme.
	
	Now, clearly, $\Isom(-)$ is a torsor of a group scheme that classify group-scheme automorphisms of $G_2$, which on its turn is a twisted form of a scheme that classify group-scheme automorphisms of $G(\tilde{\Phi},-)$; and the latter as a scheme is isomorphic to the product of $G(\tilde{\Phi},-)$ with a finite constant scheme. Thus $\Isom(-)$ is smooth and of finite type.
\end{proof}

Now consider the following subfunctor of the functor $\Isom(-)$:
$$
\Isom^\#(S)=\{\Theta\in\Isom(S)\mid \restr{(\Theta_{S}\circ (i_S\circ\ph)_*)}{E_1(R_1)}=(i_S)_*\circ\theta \}\tc
$$
where $\map{i_S}{R_2}{S}$ is a structural homomorphism

\smallskip

For every $\alpha\in\Phi$ choose a finite set $X_{\alpha}\sub G_{\alpha}^{(1)}(R_1)$ such that $X_{\alpha}$ generate $G_{\alpha}^{(1)}(R_1)$ or $G^{(1)}_{\alpha}(R_1)/G^{(1)}_{2\alpha}(R_1)$ if $\alpha$ is ultrashort) as an $R_1$-module.

{\lem\label{IsomXIsEnough} Let $S$ be an $R_2$-algebra. Let $\Theta\in\Isom(S)$ be such that for every $\alpha\in\Phi$ we have $\restr{(\Theta_{S}\circ (i_S\circ\ph)_*)}{X_{\alpha}}=\restr{((i_S)_*\circ\theta)}{X_{\alpha}}$. Then $\Theta\in\Isom^\#(S)$.}
\begin{proof}
	Clearly, for any $S$-algebra $A$ the set $(i_A\circ\ph)_*(X_{\alpha})$ generate $\leftact{\ph}{G}_{\alpha}^{(1)}(A)$ (or $\leftact{\ph}{G}^{(1)}_{\alpha}(A)/\leftact{\ph}{G}^{(1)}_{2\alpha}(A)$ if $\alpha$ is ultrashort) as an $A$-module; and $((i_A)_*\circ\theta)(X_{\alpha})$ generate $G_{\alpha}^{(2)}(A)$ or $G^{(2)}_{\alpha}(A)/G^{(2)}_{2\alpha}(A)$ if $\alpha$ is ultrashort) as an $A$-module. Thus Lemma~\ref{EgorFormulas} implies that $\Theta_A(\leftact{\ph}{G}^{(1)}_{\alpha}(A))=G^{(2)}_{\alpha}(A)$.
	
	It follows from the proof of \cite[Theorem~5]{VoronDiophantine} that there is a unique ring automorphism $\eta_A\colon A\toiso A$ such that for any $\alpha\in\Phi$, for any $x\in \leftact{\ph}{G}_{\alpha}^{(1)}(A)$ and for any $\xi\in A$ we have $\Theta_A(\xi\cdot x)=\eta_A(\xi)\cdot \Theta_A(x)$ if $\alpha$ is non-ultrashort; and $\Theta_A(\xi\cdot x)\in(\eta_A(\xi)\cdot \Theta_A(x))G_{2\alpha}^{(2)}(A)$ if $\alpha$ is ultrashort.
	
	It is easy to see that collection of all $\eta_A$ is an automorphism of the forgetful functor from the category of $S$-algebras, to the category of rings. Then by Lemma~\ref{RingStays} each $\eta_A$ is the identity automorphism. In particular, we have $\eta_{S}=\id_{S}$. Thus for any $x\in \leftact{\ph}{G}_{\alpha}^{(1)}(R_2)$ and for any $\xi\in R_2$ we have
	$$
	\Theta_{S}((i_S)_*(\xi\cdot x))=\Theta_{S}(i_S(\xi)\cdot(i_S)_*(x))=i_S(\xi)\cdot \Theta_{S}((i_S)_*(x))
	$$
	 if $\alpha$ is non-ultrashort; and
	 $$
	 \Theta_{S}((i_S)_*(\xi\cdot x))\in(i_S(\xi)\cdot \Theta_{S}((i_S)_*(x)))G_{2\alpha}^{(2)}(S)
	 $$
	  if $\alpha$ is ultrashort. In other words, for any $x\in G_{\alpha}^{(1)}(R_1)$ and for any $\xi\in R_1$ we have
	  $$
	  (\Theta_{S}\circ (i_S\circ\ph)_*)(\xi\cdot x)=(i_S\circ\ph)(\xi)\cdot (\Theta_{S}\circ (i_S\circ\ph)_*)(x)
	  $$
	  if $\alpha$ is non-ultrashort; and
	  $$
	  (\Theta_{S}\circ (i_S\circ\ph)_*)(\xi\cdot x)\in( (i_S\circ\ph)(\xi)\cdot (\Theta_{S}\circ (i_S\circ\ph)_*)(x))G_{2\alpha}^{(2)}(S)
	  $$
	   if $\alpha$ is ultrashort.
	
	Now for any $x\in X_{\alpha}$ and for any $\xi\in R_1$ we have
	\begin{multline*}
		(\Theta_{S}\circ (i_S\circ\ph)_*)(\xi\cdot x)=(i_S\circ\ph)(\xi)\cdot (\Theta_{S}\circ (i_S\circ\ph)_*)(x)=(i_S\circ\ph)(\xi)\cdot ((i_S)_*\circ\theta)(x)=\\=(i_S)_*(\ph(\xi)\cdot \theta(x))=((i_S)_*\circ\theta)(\xi\cdot x)
	\end{multline*}
	if $\alpha$ is non-ultrashort; and similarly
	$$
	(\Theta_{S}\circ (i_S\circ\ph)_*)(\xi\cdot x)\in ((i_S)_*\circ\theta)(\xi\cdot x)G_{2\alpha}^{(2)}(S)
	$$
	if $\alpha$ is ultrashort.
	
	Since addition in the modules $G_{\alpha}^{(1)}(R_1)$ and $G_{\alpha}^{(1)}(R_1)/G_{2\alpha}^{(1)}(R_1)$ and in the modules $G_{\alpha}^{(2)}(S)$ and $G_{\alpha}^{(2)}(S)/G_{2\alpha}^{(2)}(S)$ are defined by group operation, we conclude that the homomorphisms $\Theta_{S}\circ (i_S\circ\ph)_*$ and $(i_S)_*\circ\theta$ induce the same maps $G_{\alpha}^{(1)}(R_1)\to G_{\alpha}^{(2)}(S)$ and  $G_{\alpha}^{(1)}(R_1)/G_{2\alpha}^{(1)}(R_1)\to G_{\alpha}^{(2)}(S)/G_{2\alpha}^{(2)}(S)$. Therefore, by Proposition~\ref{UltrashortFix} we have $\restr{(\Theta_{S}\circ (i_S\circ\ph)_*)}{E_1(R_1)}=(i_S)_*\circ\theta$, i.e $\Theta\in\Isom^\#(S)$.
	\end{proof}

{\lem\label{IsomIsFinitePresented} The functor $\Isom^\#(-)$ is a finitely presented affine scheme over $R_2$.}
\begin{proof}
	By Lemma~\ref{IsomGood} the scheme $\Isom(-)$ is affine and smooth, so particular it is finitely presented. Now by Lemma~\ref{IsomXIsEnough} $\Isom^\#(-)$ is a closed subscheme of $\Isom(-)$ defined by finite number of equations. Hence it is affine and finitely presented.
\end{proof}

{\lem\label{Uniqueness} For any $R_2$-algebra $S$ we have $|\Isom^\#(S)|\le 1$.}
\begin{proof}
	Let $\Theta$,$\Theta'\in \Isom^\#(S)$. Let $A$ be an $S$-algebra. Similarly to how we have done it in the proof of Lemma~\ref{IsomXIsEnough}, we can use Lemma~\ref{RingStays} in order to prove that for any $x\in \leftact{\ph}{G}_{\alpha}^{(1)}(A)$ and for any $\xi\in A$ we have $\Theta_{A}(\xi\cdot x)=\xi\cdot \Theta_{A}(x)$ if $\alpha$ is non-ultrashort; and $\Theta_{A}(\xi\cdot x)\in(\xi\cdot \Theta_{A}(x))G_{2\alpha}^{(2)}(A)$ if $\alpha$ is ultrashort. The same holds for $\Theta'$.
	
	Since addition in the modules $\leftact{\ph}{G}_{\alpha}^{(1)}(A)$ and $\leftact{\ph}{G}_{\alpha}^{(1)}(A)/\leftact{\ph}{G}_{2\alpha}^{(1)}(A)$ is defined by group operation; and since those modules are generated by the images of $G_{\alpha}^{(1)}(R_1)$, it follows that the homomorphisms $\Theta$ and $\Theta'$ induce the same bijections $\leftact{\ph}{G}_{\alpha}^{(1)}(A)\toiso G_{\alpha}^{(2)}(A)$ and  $\leftact{\ph}{G}_{\alpha}^{(1)}(A)/\leftact{\ph}{G}_{2\alpha}^{(1)}(A)\toiso G_{\alpha}^{(2)}(A)/G_{2\alpha}^{(2)}(A)$. Therefore, by Proposition~\ref{UltrashortFix} we obtain that $\Theta$ and $\Theta'$ coincide on the elementary subgroup of $\leftact{\ph}{G_1}(A)$. Since every element of $\leftact{\ph}{G_1}(A)$ becomes elementary in a faithfully flat extension of $A$, it follows that $\Theta_A=\Theta'_A$. Since $A$ here is any $S$-algebra, it follows that $\Theta=\Theta'$.
\end{proof}

{\lem\label{IsomOverFields} For any maximal ideal $J\unlhd R_2$ we have $\Isom^\#(R_2/J)\ne\emp$.}

\begin{proof}
	Let $I\unlhd R_1$ be a maximal ideal that correspond to $I$ by correspondence from Lemma~\ref{MaxCorr}. Then $\theta$ induces an isomorphism $\overline{\theta}\colon E_1(R_1/I)\toiso E_2(R_2/J)$. By Proposition~\ref{BT} $\overline{\theta}$ comes from a field isomorphism $\overline{\ph}\colon R_1/I\toiso R_2/J$ and an $R_2/J$-group-scheme isomorphism $\Theta\colon \leftact{\overline{\ph}}{(G_1)_{R_1/I}}\toiso (G_2)_{R_2/J}$ (all the exceptional cases are excluded by assumptions of Theorem~\ref{main} and the assumption $\tilde{\Phi_1}=\tilde{\Phi_2}$).
	
	In order to show that $\Theta$ can be viewed as an element of $\Isom^\#(R_2/J)$ we only must prove that $\overline{\ph}\circ \rho_I=\rho_J\circ \ph$, where $\map{\rho_I}{R_1}{R_1/I}$ and $\map{\rho_J}{R_2}{R_2/J}$ are the projections on the quotients.
	
	Similarly to how we have done it in the proof of Lemma~\ref{IsomXIsEnough}, we can use Lemma~\ref{RingStays} in order to prove that for any $\overline{x}\in G_{\alpha}^{(1)}(R_1/I)$ and for any $\overline{\xi}\in R_1/I$ we have
	$$
	(\Theta_{R_2/J}\circ\overline{\ph}_*)(\overline{\xi}\cdot \overline{x})=\overline{\ph}(\overline{\xi})\cdot (\Theta_{R_2/J}\circ\overline{\ph}_*)(\overline{x})
	$$
	 if $\alpha$ is non-ultrashort.
	
	Therefore, for any $x\in G_{\alpha}^{(1)}(R_1)$ and any $\xi\in R_1$ we have
	\begin{multline*}
		\overline{\ph}(\rho_I(\xi))\cdot (\Theta_{R_2/J}\circ\overline{\ph}_*)((\rho_I)_*(x))=	(\Theta_{R_2/J}\circ\overline{\ph}_*)(\rho_I(\xi)\cdot (\rho_I)_*(x))=\overline{\theta}(\rho_I(\xi)\cdot (\rho_I)_*(x))=\\=\overline{\theta}((\rho_I)_*(\xi\cdot x))=(\rho_J)_*(\theta(\xi\cdot x))=(\rho_J)_*(\ph(\xi)\cdot\theta(x))=\rho_J(\ph(\xi))\cdot (\rho_J)_*(\theta(x))=\\=\rho_J(\ph(\xi))\cdot\overline{\theta}((\rho_I)_*(x))=\rho_J(\ph(\xi))\cdot (\Theta_{R_2/J}\circ\overline{\ph}_*)((\rho_I)_*(x))\tp
	\end{multline*}
	
	Since the map $\map{\Theta_{R_2/J}\circ\overline{\ph}_*}{G_{\alpha}^{(1)}(R_1/I)}{ G_{\alpha}^{(2)}(R_2/J)}$ is an isomorphism and the map $\map{(\rho_I)_*}{G_{\alpha}^{(1)}(R_1)}{G_{\alpha}^{(1)}(R_1/I)}$ is a surjection, it follows that for any $\xi\in R_1$ and any $y\in G_{\alpha}^{(2)}(R_2/J)$ we have
	$$
	\overline{\ph}(\rho_I(\xi))\cdot y=\rho_J(\ph(\xi))\cdot y\tp
	$$
	Since $G_{\alpha}^{(2)}(R_2/J)$ is a non-trivial vector space, it follows that $\overline{\ph}(\rho_I(\xi))=\rho_J(\ph(\xi))$.
\end{proof}

{\lem\label{IsomSmooth} The scheme $\Isom^\#(-)$ is smooth over $R_2$.}
\begin{proof}
	Let $S_0$ be a faithfully flat finitely presented $R_2$-algebra, such that the isotropic pinning $(T^{(2)},\Phi)$ on $G_2$ after a base change to $S_0$ becomes contained in a pinning. Since smoothness is a fppf-local property, it is enough to show that $\Isom^\#(-)_{S_0}$ is smooth over~$S_0$.
	
	By Lemma~\ref{IsomIsFinitePresented} the scheme $\Isom^\#(-)$ is locally finitely presented; hence, so is $\Isom^\#(-)_{S_0}$. Thus by \cite[Ch I, \S 4, Item 4.5]{DemazureGabriel} it is enough to show that for any $S_0$-algebra $S$ and any ideal $I\unlhd S$ such that $I^2=0$ the reduction map $\Isom^\#(S)\to\Isom^\#(S/I)$ is surjective.
	
	Now let $\overline{\Theta}\in \Isom^\#(S/I)$. Since by Lemma~\ref{IsomGood} the scheme $\Isom(-)$ is smooth, there is $\Theta'\in \Isom(S)$ such that its reduction modulo $I$ is $\overline{\Theta}$. Using Proposition~\ref{InfinitesiamiPartI} for the map $(\Theta'_S\circ (i_S\circ\ph)_*)$ (here all the umbrella assumptions of Section~\ref{InfinitesimalSecI} are satisfied because $\overline{\Theta}\in  \Isom^\#(S/I)$: if follows from Lemma~\ref{EgorFormulas} that $\overline{\Theta}_{S/I}$ preserves root subgroups; and since elements $h_{\alpha}$ are defined through the root subgroups $\overline{\Theta}_{S/I}$ must preserve them as well) we conclude that there exists $g\in G_2(S,I)\simeq G(\tilde{\Phi},S,I)$ such that $(\Theta'_S\circ (i_S\circ\ph)_*)^g$ preserves root subgroups.
	
	Therefore, we can assume without loss of generality that $\Theta'_S\circ (i_S\circ\ph)_*$ preserves root subgroups, which by Lemma~\ref{EgorFormulas} implies that $\Theta'_S(\leftact{\ph}{G}_{\alpha}^{(1)}(S))=G_{\alpha}^{(2)}(S)$ for all $\alpha\in\Phi$.
	
	 Denote by $\map{(\Theta'_S)_{\alpha}}{\leftact{\ph}{G}_{\alpha}^{(1)}(S)}{G_{\alpha}^{(2)}(S)}$ (or $\map{(\Theta'_S)_{\alpha}}{\leftact{\ph}{G}_{\alpha}^{(1)}(S)/\leftact{\ph}{G}_{2\alpha}^{(1)}(S)}{G_{\alpha}^{(2)}(S)/G_{2\alpha}^{(2)}(S)}$ if $\alpha$ is ultrashort) the group isomorphisms induced by $\Theta'_S$ for every $\alpha\in\Phi$. Similarly to how we have done it in the proof of Lemma~\ref{IsomXIsEnough}, we can use Lemma~\ref{RingStays} in order to prove that $(\Theta'_S)_{\alpha}$ are in fact isomorphisms of $S$-modules.
	
	Now let $\map{\theta_{\alpha}}{G_{\alpha}^{(1)}(R_1)}{G_{\alpha}^{(2)}(R_2)}$ (or $\map{\theta_{\alpha}}{G_{\alpha}^{(1)}(R_1)/G_{2\alpha}^{(1)}(R_1)}{G_{\alpha}^{(2)}(R_2)/G_{2\alpha}^{(2)}(R_2)}$ if $\alpha$ is ultrashort) be the group isomorphisms induced by $\theta$. Then we have group isomorphisms $\map{\theta_{\alpha}\circ\ph_*^{-1}}{\leftact{\ph}{G}_{\alpha}^{(1)}(R_2)}{G_{\alpha}^{(2)}(R_2)}$ (or $\map{\theta_{\alpha}\circ\ph_*^{-1}}{\leftact{\ph}{G}_{\alpha}^{(1)}(R_2)/\leftact{\ph}{G}_{2\alpha}^{(1)}(R_2)}{G_{\alpha}^{(2)}(R_2)/G_{2\alpha}^{(2)}(R_2)}$ if $\alpha$ is ultrashort). It follow from the definition of $\ph$ that the last ones are the homomorphisms of $R_2$-modules.
	
	Since we have a natural identification $\leftact{\ph}{G}_{\alpha}^{(1)}(S)=\leftact{\ph}{G}_{\alpha}^{(1)}(R_2)\otimes S$ (or $\leftact{\ph}{G}_{\alpha}^{(1)}(S)/\leftact{\ph}{G}_{2\alpha}^{(1)}(S)=\leftact{\ph}{G}_{\alpha}^{(1)}(R_2)/\leftact{\ph}{G}_{\alpha}^{(1)}(R_2)\otimes S$ if $\alpha$ is ultrashort); and $G_{\alpha}^{(2)}(S)=G_{\alpha}^{(2)}(R_2)\otimes S$ (or $G_{\alpha}^{(2)}(S)/G_{2\alpha}^{(2)}(S)=G_{\alpha}^{(2)}(R_2)/G_{\alpha}^{(2)}(R_2)\otimes S$ if $\alpha$ is ultrashort), it follows that the isomorphisms $\theta_{\alpha}\circ\ph_*^{-1}$ induce the isomorphisms $\map{(\theta_{\alpha}\circ\ph_*^{-1})_S}{\leftact{\ph}{G}_{\alpha}^{(1)}(S)}{G_{\alpha}^{(2)}(S)}$ (or $\map{(\theta_{\alpha}\circ\ph_*^{-1})_S}{\leftact{\ph}{G}_{\alpha}^{(1)}(S)/\leftact{\ph}{G}_{2\alpha}^{(1)}(S)}{G_{\alpha}^{(2)}(S)/G_{2\alpha}^{(2)}(S)}$ if $\alpha$ is ultrashort).
	
	Now set $\sigma_{\alpha}=(\theta_{\alpha}\circ\ph_*^{-1})_S\circ (\Theta'_S)_{\alpha}^{-1}$. So $\sigma_{\alpha}$ is an automorphism of $G_{\alpha}^{(2)}(S)$ (or $G_{\alpha}^{(2)}(S)/G_{2\alpha}^{(2)}(S)$ if $\alpha$ is ultrashort) as an $S$-module.
	
	\medskip
	
	\underline{Sublemma} The collection $(\sigma_{\alpha})_{\alpha\in\Phi}$ defines via isomorphisms $t_{\alpha}\colon\g_{\alpha}(S)\toiso G_{\alpha}^{(2)}(S)$ (or $t_{\alpha}\colon \g_{\alpha}(S)\toiso G_{\alpha}^{(2)}(S)/G_{2\alpha}^{(2)}(S)$ if $\alpha$ is ultrashort) an element in $\Diag(S)$, where $\Diag(-)$ is a scheme from Lemma~\ref{Clopen} and $t_{\alpha}$ are the exponent maps, which are known to be well defined.
	
	\medskip
	
	\underline{Proof of Sublemma.} Let $\alpha$,$\beta\in\Phi$ with $\alpha+\beta\in\Phi$ are either independent or equal. For any $x\in \g_{\alpha}(S)$ and $y\in\g_{\beta}(S)$ we have
	$$
	[t_{\alpha}(x),t_{\beta}(y)]\in t_{\alpha+\beta}([x,y])\prod_{i+j\ge 3}G_{i\alpha+j\beta}^{(2)}(S)\text{;}
	$$
	hence, we must prove that for any $x\in \g_{\alpha}(S)$ and $y\in\g_{\beta}(S)$ we have
	$$
	[\sigma_{\alpha}(t_{\alpha}(x)),\sigma_{\beta}(t_{\beta}(y))]\in \sigma_{\alpha+\beta}(t_{\alpha+\beta}([x,y]))\prod_{i+j\ge 3}G_{i\alpha+j\beta}^{(2)}(S)
	$$
	(when somebody is ultrashort replace ''$\in$'' by ''$\sub$'').
	
	We may assume that $x\in (\Theta'_S)_{\alpha}(\leftact{\ph}{G}_{\alpha}^{(1)}(R_2))$ (resp. $x\in (\Theta'_S)_{\alpha}(\leftact{\ph}{G}_{\alpha}^{(1)}(R_2)/\leftact{\ph}{G}_{2\alpha}^{(1)}(R_2))$) and $y\in (\Theta'_S)_{\beta}(\leftact{\ph}{G}_{\beta}^{(1)}(R_2))$ (resp. $x\in (\Theta'_S)_{\beta}(\leftact{\ph}{G}_{\beta}^{(1)}(R_2)/\leftact{\ph}{G}_{2\beta}^{(1)}(R_2))$) because these subsets generate $G_{\alpha}^{(2)}(S)$ and $G_{\beta}^{(2)}(S)$ (resp. $G_{\alpha}^{(2)}(S)/G_{2\alpha}^{(2)}(S)$ and $G_{\beta}^{(2)}(S)/G_{2\beta}^{(2)}(S)$) as $S$-modules. In this case the statement follows from the fact that $\Theta'_{S}$ and $\theta\circ\ph_*^{-1}$ are group homomorphisms.
	
	\underline{Sublemma is proved.}
	
	\medskip
	
	Now the fact that $\overline{\Theta}\in \Isom^\#(S/I)$ implies that each $\sigma_{\alpha}$ is congruent to the identity operator modulo $I$; hence by Lemma~\ref{Clopen} the collection $(\sigma_{\alpha})_{\alpha\in\Phi}$  must come from an adjoint action of the element $g\in L(S,I)\le G(\tilde{\Phi},S,I)=G_2(S,I)$.
	
	Now set $\Theta=(\Theta')^g$. It is easy to see that the maps $\Theta_{S}\circ (i_S\circ\ph)_*$ and $(i_S)_*\circ\theta$ induce the same maps $G_{\alpha}^{(1)}(R_1)\to G_{\alpha}^{(2)}(S)$ and  $G_{\alpha}^{(1)}(R_1)/G_{2\alpha}^{(1)}(R_1)\to G_{\alpha}^{(2)}(S)/G_{2\alpha}^{(2)}(S)$. Therefore, by Proposition~\ref{UltrashortFix} we have $\restr{(\Theta_{S}\circ (i_S\circ\ph)_*)}{E_1(R_1)}=(i_S)_*\circ\theta$; hence, we have $\Theta\in\Isom^\#(S)$; and, clearly, the reduction of $\Theta$ modulo $I$ is $\overline{\Theta}$.
	
\end{proof}

\section{The proof of item (2) of Theorem~\ref{main}}
\label{ProofFinaleSec}

The morphism $\Isom^\#(-)\to \Spec R_2$ if flat by Lemma~\ref{IsomSmooth}; hence, the map on the underlying topological spaces has an open image. By Lemma~\ref{IsomOverFields} this image contains all the closed points. Since $\Spec R_2$ is affine this implies that the map on the underlying topological spaces is surjective. Therefore, the ring $R_2[\Isom^\#]$ of regular functions on $\Isom^\#(-)$ is a faithfully flat $R_2$-algebra.

This implies that the structural homomorphism $R_2\to R_2[\Isom^\#]$ is the equalizer of two homomorphisms from $R_2[\Isom^\#]$ to $R_2[\Isom^\#]\otimes_{R_2}R_2[\Isom^\#]$. These homomorphisms correspond to the elements of $\Isom^\#(R_2[\Isom^\#]\otimes_{R_2}R_2[\Isom^\#])$; hence, by Lemma~\ref{Uniqueness} they are equal; hence the structural homomorphism $R_2\to R_2[\Isom^\#]$ is an isomorphism. The inverse homomorphism correspond to an element of $\Isom^\#(R_2)$. Therefore, we have $\Isom^\#(R_2)\ne\emp$, which by the definition of $\Isom^\#(-)$ finishes the proof.

\section{Generalization to bigger subgroups}
\label{CorollarySec}

The lemma below does not require the group scheme to have an isotropic pinning. So this is in greater generality than the rest of the paper. Here we use the notion of the elementary subgroup $E(R)$ from~\cite{PetStavIso} and the notation $E(R,I)$ stands for the relative elementary subgroup as defined in \cite{StavStep}.

{\lem\label{CharacteristicElementary} Let $G$ be an absolutely simple adjoint group scheme with a common root datum of the geometric fibers over the ring $R$. Assume that $G$ has isotropic rank at least 2. If the absolute root system of $G$ is doubly laced, assume that $2\in R^*$; and if the absolute root system of $G$ is of type $G_2$, assume that $6\in R^*$. Suppose that we have an intermediate subgroup $E(R)\le H\le G(R)$. Then $E(R)$ is the smallest by inclusion among all the subgroups $K\le H$ that satisfy the following properties:

1) $K$ is a normal closure in $H$ of a finitely generated subgroup;

2) $K=[K,K]$;

3) the centralizer of $K$ in $H$ is trivial.}

\begin{proof}
First let us show that $E(R)$ satisfies the properties (1)--(3). By results of~\cite{PetStavIso} the subgroup $E(R)$ is normal in $G(R)$. Now let $\alpha$ be a relative root, and consider the closed embedding $\map{X_{\alpha}}{W(V_{\alpha})}{G}$ from~\cite[Theorem~2]{PetStavIso}, where $V_{\alpha}$ is a finitely generated projective module of constant non-zero rank; and $W(V_{\alpha})$ is the $R$-scheme defined by $W(V_{\alpha})(S)=V_{\alpha}\otimes_R S$. Let $v_1$,$\ldots$,$v_n\in V_{\alpha}=W(V_{\alpha})(R)$ be the generators of $V_{\alpha}$. We claim that $E(R)$ is s a normal closure in $H$ of a subgroup generated by $X_{\alpha}(v_1)$,$\ldots$,$X_{\alpha}(v_n)$. Indeed, since the image of $X_{\alpha}$ is the root subgroup, we have $X_{\alpha}(v_1)$,$\ldots$,$X_{\alpha}(v_n)\in E(R)$. Now suppose that $K\unlhd H$ is a normal subgroup that contains $X_{\alpha}(v_1)$,$\ldots$,$X_{\alpha}(v_n)$, but does not contain $E(R)$. By~\cite[Theorem 1.1]{StavStep} we have $K\le G(R,I)$ for some proper ideal $I$. Thus, we have $v_1$,$\ldots$,$v_n\in IV_{\alpha}$, which contradict the fact that they generate $V_{\alpha}$. Therefore, $E(R)$ satisfies~(1). Further $E(R)$ satisfies~(2) by~\cite{LuzStav}; and it satisfies~(3) by~\cite{KulikStav}.
	
	Now let $K\le H$ be a subgroup that satisfies (1)--(3); we must prove that $E(R)\le K$. Since $K$ is normal in $H$, it follows by~\cite[Theorem 1.1]{StavStep} that $E(R,I)\le K\le G(R,I)$ for some ideal $I\unlhd R$. Assume that $K$ is a normal closure in $H$ of a  subgroup generated by elements $g_1$,$\ldots$,$g_n$. Let $f_1$,$\ldots$,$f_m$ be the generators of the algebra $R[G]$ of regular functions on $G$. Clearly, $I$ is the smallest ideal such that $K\le G(R,I)$; hence $I$ must be generated by $f_i(g_j)-a(f_i)$ for all $1\le i\le m$ and $1\le j\le n$, where $\map{a}{R[G]}{R}$ is the augmentation homomorphism. Thus the ideal $I$ is finitely generated.
	
	By (2) we have $E(R,I)\le K=[K,K]\le[G(R,I),G(R,I)]\le G(R,I^2)$; hence, we have $I=I^2$. Since $I$ is finitely generated, it follows by Nakayama's lemma that $I$ is generated by some idempotent $e\in R$. Thus $K$ commutes with $E(R,(1-e))$; hence, (3) implies that $e=1$; hence, we have $I=R$; hence, we have $E(R)\le K$.
\end{proof}

{\cor\label{Bigger} Under the assumptions of Theorem~\ref{main} suppose that we have intermediate subgroups $E_1(R_1)\le H_1\le G_1(R_1)$ and $E_2(R_2)\le H_2\le G_2(R_2)$. Let $\theta\colon H_1\toiso H_2$ be the isomorphism of abstract groups.
	
	Then
	
	\begin{enumerate}
		\item If $\tilde{\Phi_1}$ is not isomorphic to $\tilde{\Phi_2}$, then $\tilde{\Phi_1}=A_3$, $\tilde{\Phi_2}=B_2$, $R_1/\M\simeq \F_2$ for all maximal ideals $\M\unlhd R_1$ and $R_2/\M\simeq \F_3$ for all maximal ideals $\M\unlhd R_2$.
		
		\item  If $\tilde{\Phi_1}=\tilde{\Phi_2}$, then there exists a ring isomorphism $\ph\colon R_1\toiso R_2$ and an $R_2$-group-scheme isomorphism $\Theta\colon \phan^\ph G_1\toiso G_2$ such that $\theta=\restr{(\Theta_{R_2}\circ \ph_*)}{H_1}$.
	\end{enumerate}
}
\begin{proof}
	Since Lemma~\ref{CharacteristicElementary} characterize the elementary subgroup in purely group theoretic terms, it follows that $\theta(E_1(R_1))=E_2(R_2)$. Item (1) thus follows from the item (1) of Theorem~\ref{main}.
	
	Now suppose that $\tilde{\Phi_1}=\tilde{\Phi_2}$. By Theorem~\ref{main} there exists a ring isomorphism $\ph\colon R_1\toiso R_2$ and an $R_2$-group-scheme isomorphism $\Theta\colon \phan^\ph G_1\toiso G_2$ such that $\restr{\theta}{E_1(R_1)}=\restr{(\Theta_{R_2}\circ \ph_*)}{E_1(R_1)}$. Let us prove that actually we have $\theta=\restr{(\Theta_{R_2}\circ \ph_*)}{H_1}$.
	
	Let $h\in H_1$ and $g\in E_1(R_1)$. Using the fact that $\restr{\theta}{E_1(R_1)}=\restr{(\Theta_{R_2}\circ \ph_*)}{E_1(R_1)}$ and that $E_1(R)$ is normal in $H_1$ we obtain
	\begin{multline*}
	(\Theta_{R_2}\circ \ph_*)(h)(\theta(h))^{-1}\theta(g)\theta(h)((\Theta_{R_2}\circ \ph_*)(h)^{-1})=(\Theta_{R_2}\circ \ph_*)(h)\theta(h^{-1}gh)((\Theta_{R_2}\circ \ph_*)(h)^{-1})=\\=(\Theta_{R_2}\circ \ph_*)(h)(\Theta_{R_2}\circ \ph_*)(h^{-1}gh)((\Theta_{R_2}\circ \ph_*)(h)^{-1})=(\Theta_{R_2}\circ \ph_*)(hh^{-1}ghh^{-1})=\\=(\Theta_{R_2}\circ \ph_*)(g)=\theta(g)\tp
\end{multline*}
	
	Therefore, $(\Theta_{R_2}\circ \ph_*)(h)(\theta(h))^{-1}$ commutes with $\theta(g)$. Since $\theta(E_1(R_1))=E_2(R_2)$, it follows that $(\Theta_{R_2}\circ \ph_*)(h)(\theta(h))^{-1}$ commutes with $E_2(R_2)$; hence, by~\cite{KulikStav} we conclude that $(\Theta_{R_2}\circ \ph_*)(h)(\theta(h))^{-1}=e$; hence, $(\Theta_{R_2}\circ \ph_*)(h)=\theta(h)$.
\end{proof}

\section{Explanation of the assumption ($d$)}
\label{AssumptiondSec}

Recall that the assumption (d) in the second bullet in the statement of Theorem~\ref{main} reads: ''if $\Phi$ is of type $C$ or $BC$ and the map $u$ is not a bijection, then for every pair of orthogonal short roots $\beta$,$\beta'\in\Phi$ with their sum being a long root the corresponding Weyl elements $w_{\beta}=a_{\beta}b_{\beta}c_{\beta}$ and $w_{\beta'}=a_{\beta'}b_{\beta'}c_{\beta'}$ can be chosen so that $b_{\beta}$ commutes with $b_{\beta'}$''. Since the other assumptions require the map $u$ to come from a Tits index and the system $\Phi$ to have rank at least 2, it follows that the assumption (d) is only relevant if the Tits index in question is one of the following: $\leftact{2\!\!}{A}_{n,r}^{(d)}$; $B_{n,2}$; $C_{n,r}^{(d)}\, (d>1)$; $\leftact{1\!\!}{D}_{n,r}^{(d)}\,(d>1)$; $\leftact{1\!\!}{D}_{n,2}^{(1)}$; $\leftact{2\!\!}{D}_{n,r}^{(d)}\,(d>1)$; $\leftact{2\!\!}{D}_{n,2}^{(1)}$; $\leftact{2\!\!}{E}_{6,2}^{16'}$; $E_{7,2}^{31}$; $E_{8,2}^{66}$ (in each case $n\ge 3$ and $r\ge 2$ is assumed; and $n\ge 4$ is assumed for type $D$). In this Section we explain the meaning of the assumption~(d) in each of these cases.

{\prop Suppose that the Tits index in question is $C_{n,r}^{(d)}\, (d>1)$; i.e. for some Azumaya algebra $A$ over $R_1$ with an $R_1$-linear orthogonal involution $\xi\mapsto\overline{\xi}$ and a non-degenerate anti-hermitian form over $A$ of rank $\tfrac{2n}{d}$ and Witt index $r$ the isotropic pinning in question comes from an isomorphism $G_1(-)\simeq \mathrm{PU}(h\otimes_{R_1}-)$, where $\mathrm{PU}$ means the projective unitary group. Then the assumption (d) is equivalent to the assumption ''the algebra $A$ contains an invertible anti-hermitian element''.}

\begin{proof}
	In this case the short root subgroups in $G_1(R_1)$ are naturally identified with the additive group of $A$; and for a pair of orthogonal short roots $\beta$,$\beta'\in\Phi$ with their sum being a long root the commutator between the corresponding subgroups correspond to the operation $\xi\overline{\zeta}+\zeta\overline{\xi}$ for $\xi$,$\zeta\in A$. The element from $A$ can correspond to the middle component of a Weyl element iff it is invertible; hence the assumption (d) is equivalent to ''there are $\xi$,$\zeta\in A^*$ such that $\xi\overline{\zeta}+\zeta\overline{\xi}=0$''. In that case $\xi\overline{\zeta}$ is an invertible anti-hermitian element. Conversely, if $\xi\in A^*$ is anti-hermitian, then we can take $\zeta=1$.
\end{proof}

{\prop Suppose that the Tits index in question is  one of the following: $\leftact{2\!\!}{A}_{n,r}^{(d)}$,  $\leftact{1\!\!}{D}_{n,r}^{(d)}\,(d>1)$, or $\leftact{2\!\!}{D}_{n,r}^{(d)}\,(d>1)$. Then the assumption (d) is satisfied automatically.}

\begin{proof}
	The proof goes the same as for the previous proposition, except this time the commutator correspond to the operation $\xi\overline{\zeta}-\zeta\overline{\xi}$ for the elements of the corresponding Azumaya algebra; hence, we can take $\xi=\zeta=1$.
\end{proof}

Denote by $\mathcal{H}$ the hyperbolic plane over $R_1$, i.e. the module $R_1^2$ with the quadratic form given by multiplication of the coordinates.

{\prop Suppose that the Tits index in question is one of the following: $B_{n,2}$, $\leftact{1\!\!}{D}_{n,2}^{(1)}$, or $\leftact{2\!\!}{D}_{n,2}^{(1)}$; i.e for some $R_1$ module $V$ with a semi-regular quadratic form $q$ the isotropic pinning in question comes from an isomorphism $G_1(-)\simeq \mathrm{PSO}((\mathcal{H}\oplus\mathcal{H}\oplus V)\otimes_{R_1}-)$, where $\mathrm{PSO}$ means the projective special orthogonal group; and the direct sum $\mathcal{H}\oplus\mathcal{H}\oplus V$ is the orthogonal sum. Then the assumption (d) is equivalent to the assumption: ''there are orthogonal elements $v,w\in V$ with $q(v)$,$q(w)\in R_1^*$'' {\rm(}recall that the orthogonality means that $q(v+w)=q(v)=q(w)${\rm)}.}

\begin{proof}
	The proof is straightforward: the short root subgroups in $G_1(R_1)$ are naturally identified with the additive group of $V$; for a pair of orthogonal short roots $\beta$,$\beta'\in\Phi$ with their sum being a long root the commutator between the corresponding subgroups correspond to the operation $q(v+w)-q(v)-q(w)$; and the element $v\in V$ can correspond to the middle component of a Weyl element iff $q(v)\in R^*$ (see~\cite[Section 7.3B]{WiedemannRootGraded}).
\end{proof}

We will now explain how to decipher the assumption (d), if the Tits index in question is $\leftact{2\!\!}{E}_{6,2}^{16'}$, $E_{7,2}^{31}$ or $E_{8,2}^{66}$. In each of these cases we have $\Phi=BC_2$. Now let $\Phi_{\mathrm{nus}}\le\Phi$ be the subsystem of non-ultrashort roots, and let $\tilde{\Phi}_{\mathrm{nus}}=u^{-1}(\Phi_{\mathrm{nus}}\cup\{0\})\cap\tilde{\Phi}$. Then the subsystem subgroup for the $\tilde{\Phi}_{\mathrm{nus}}$ is defined over $R_1$ and it comes with an isotropic pinning with the root system $\Phi_{\mathrm{nus}}$. The corresponding map $\tilde{\Phi}_{\mathrm{nus}}\cup\{0\}\to \Phi_{\mathrm{nus}}\cup\{0\}$ is the restriction of $u$; and it is easy to see that it comes from the Tits index $\leftact{2\!\!}{D}_{5,2}^{(1)}$, $\leftact{1\!\!}{D}_{6,2}^{(1)}+\leftact{1\!\!}{A}_{1,0}^{(2)}$, or $\leftact{1\!\!}{D}_{8,2}^{(1)}$. In the second case the subsystem subgroup for $D_6\le E_7$ is defined over $R_1$, because we can define it as an fppf-sheafification of the subgroup generated by non-ultrashort root elements. Either way since the assumption $(d)$ depends only on non-ultrashort root subgroups, deciphering it for the Tits indexes $\leftact{2\!\!}{E}_{6,2}^{16'}$, $E_{7,2}^{31}$ and $E_{8,2}^{66}$ reduces to deciphering it for the Tits indexes $\leftact{1\!\!}{D}_{6,2}^{(1)}$, and $\leftact{1\!\!}{D}_{8,2}^{(1)}$, which is done in the previous proposition.


\end{document}